\documentclass[11pt]{article}
\usepackage{times, amsfonts, amsmath, amssymb, latexsym, setspace, epsfig}
\usepackage{subfigure, graphics}
\graphicspath{{figs/}}
\usepackage{multirow}
\usepackage{algorithm,algorithmic}
\usepackage{array}
\usepackage{authblk}
\newcolumntype{P}[1]{>{\centering\arraybackslash}p{#1}}

\def\bs{\boldsymbol}

\def\Pr{{\rm I\!P}}
\def\be{\begin{equation}}
\def\ee{\end{equation}}
\def\bea{\begin{eqnarray*}}
\def\eea{\end{eqnarray*}}
\def\bean{\begin{eqnarray}}
\def\eean{\end{eqnarray}}
\def\nn{\nonumber}
\def\nin{\noindent}

\def\ra{\rightarrow}

\def\Bl{\Bigl}
\def\Br{\Bigr}

\def\R{{\bf R}}

\def\alp{\alpha}

\def\AA{{\cal A}}
\def\CC{{\cal C}}
\def\II{{\cal I}}
\def\KK{{\cal K}}
\def\x{{\bs x}}
\def\g{{\bs g}}
\def\bsell{{\bs \ell}}

\newtheorem{Theorem}{Theorem}
\newtheorem{Proposition}{Proposition}

\newtheorem{Lemma}{Lemma}
\newtheorem{Fact}{Fact}

\begin{document}
\setstretch{1.1}

\title{Confidence bands for a log-concave density}

\author[1]{Guenther Walther\thanks{Research supported by NSF grants DMS-1501767 and DMS-1916074}}
\author[1,2]{Alnur Ali\thanks{Research supported by the Intelligence Community Postdoctoral Research Fellowship Program}}
\author[2,3]{Xinyue Shen\thanks{Research supported by the National Key R\&D Program of China with grant No. 2018YFB1800800,
the Key Area R\&D Program of Guangdong Province with grant No. 2018B030338001,
the Shenzhen Outstanding Talents Training Fund,
and the Guangdong Research Project No. 2017ZT07X152}}
\author[2]{Stephen Boyd} 
\affil[1]{Department of Statistics, Stanford University}
\affil[2]{Department of Electrical Engineering, Stanford University}
\affil[3]{Shenzhen Research Institute of Big Data and Future Network of Intelligence Institute, Chinese
University of Hong Kong \authorcr
\texttt{\{gwalther,alnurali,xinyues,boyd\}@stanford.edu}}

\date{revised May 2022}

\maketitle

\begin{abstract}
We present a new approach for inference about a univariate log-concave distribution: Instead of using
the method of maximum likelihood, we propose to incorporate the log-concavity constraint in
an appropriate nonparametric confidence set for the cdf $F$.
This approach has the advantage that it automatically provides a measure of statistical uncertainty and it thus overcomes
a marked limitation of the maximum likelihood estimate. In particular, we show how to construct 
confidence bands for the density that have a finite sample guaranteed confidence level.
The nonparametric confidence set for $F$ which we introduce here has attractive computational
and statistical properties: It allows to bring modern tools from optimization to bear on
this problem via difference of convex programming, and it results in optimal
statistical inference. We show that the width of the resulting confidence bands
converges at nearly the parametric $n^{-\frac{1}{2}}$ rate when the log density is $k$-affine.
\end{abstract}

\section{Introduction}

Statistical inference under shape constraints has been the subject of continued considerable
research activity. Imposing shape constraints on the inference about a function $f$, that is,
assuming that the function satisfies certain qualitative properties, such as monotonicity or
convexity on certain subsets of its domain, has two main motivations: First, such shape
constraints may directly arise from the problem under investigation and it is then desirable
that the result of the inference reflects this fact. The second reason is that alternative
nonparametric estimators typically involve the choice of a tuning parameter, such as the
bandwidth in the case of a kernel estimator. A good choice for such a tuning parameter is
usually far from trivial. More importantly, selecting a tuning parameter injects a certain amount of subjectivity
into the estimator, and the resulting choice may prove quite consequential for relevant aspects
of the inference. In contrast, imposing shape constraints often allows to derive an
estimator that does not depend on a tuning parameter while at the same time exhibiting
a good statistical performance, such as achieving optimal minimax rates of convergence
to the underlying function $f$. 

This paper is concerned with inference about a univariate log-concave density, i.e. a
density of the form
$$
f(x) = \exp \phi(x),
$$
where $\phi: \R \ra [–\infty,\infty)$ is a concave function.
It was argued in Walther~(2002) that log-concave densities
represent an attractive and useful nonparametric surrogate for the class of Gaussian
distributions for a range of problems in inference and modeling. The appealing properties 
of this class are that it contains most commonly encountered parametric families of
unimodal densities with exponentially decaying tails and that the class is closed
under convolution, affine transformations, 
convergence in distribution and marginalization.
 In a similar vein as a normal distribution 
can be seen as a prototypical model for a homogenous population, one can use the
class of log-concave densities as a more flexible nonparametric model for this task,
from which heterogenous models can then be built e.g. via location mixtures. Historically such
homogenous distributions have been modeled with unimodal densities, but it is known that
the maximum likelihood estimate (MLE) of a unimodal density does not exist, see e.g.
Birg\'{e}~(1997). In contrast, it was shown in Walther~(2002) that the MLE of a log-concave
density exists and can be computed with readily available algorithms. 
Therefore, the class of log-concave densities is a sufficiently rich
nonparametric model while
at the same time it is small enough to allow nonparametric inference without a tuning parameter.

Due to these attractive properties there has been a considerable research activity in the last
15 years about
the statistical properties of the MLE, computational aspects, applications in modeling and
inference, and the multivariate setting. Many of the key properties of the MLE are now well
understood: Existence of the MLE was shown in Walther~(2002), consistency was proved by
Pal, Woodroofe and Meyer~(2007), and rates of convergence in certain uniform metrics
was established by D\"{u}mbgen and Rufibach~(2009). Balabdaoui, Rufibach and Wellner~(2009)
provided pointwise limit distributribution theory, while Doss and Wellner~(2016) and
Kim and Samworth~(2016) gave 
rates of convergence in the Hellinger metric, and Kim, Guntuboyina and Samworth~(2018)
proved adaptation properties. Accompanying results for the multivariate case are given
in e.g. Cule, Samworth and Stewart~(2010), Schuhmacher and D\"{u}mbgen~(2010),
Seregin and Wellner~(2010), Kim and Samworth~(2016), and Feng, Guntuboyina, Kim and Samworth~(2020).

Computation of the univariate MLE was initially approached with the Iterative Convex Minorant Algorithm, see
Walther~(2002), Pal, Woodroofe and Meyer~(2007) and Rufibach~(2007), but it appears that the
fastest algorithms currently available are the
active set algorithm given in D\"{u}mbgen, H\"{u}sler and Rufibach~(2007)
and the constrained Newton method proposed by Liu and Wang~(2018).

%Computing the multivariate
%MLE requires a different approach and only slower algorithms are available.
%Cule, Samworth and Stewart~(2010) formulated this problem as a nondifferentiable convex
%optimization problem and used Shor's $r$-algorithm to solve it. Unfortunately, the resulting
%long computation time make this approach not practical even for moderate sample sizes in
%low dimensions, which has motivated fundamental research to find efficient algorithms for
%this problem, see e.g. Axelrod et al.~(2019).
%It appears that the algorithm described in Rathke and Schn\"{o}rr~(2019) is currently the 
%most efficient approach that also comes with practically
%usable code. It computes an approximation
%to the MLE rather than the MLE itself.

Overviews of some of these results and other work involving modeling and inference with log-concave
distributions are given in the review
papers of Walther~(2009), Saumard and Wellner~(2014) and Samworth~(2018). 

Notably, the existing methodology
for estimating a log-concave density appears to be exclusively focused on the method of maximum
likelihood. Here we will employ a different methodology: We will derive a confidence band
by intersecting the log-concavity constraint with a goodness-of-fit test. One important advantage of
this approach is that such a confidence band satisfies a key principle of statistical inference:
an estimate needs to be accompanied by some measure
of standard error in order to be useful for inference. There appears to be no known method for obtaining such a confidence band via
the maximum likelihood approach. Balabdaoui et al.~(2009) construct {\em pointwise} confidence intervals for a log-concave density
based on asymptotic limit theory, which requires to estimate the second derivative of $\log f(x)$. 
Azadbakhsh et al.~(2014) compare several methods for estimating this nuisance parameter and they report that this task is 
quite difficult. An alternative approach is given by Deng et al.~(2020).
In Section~\ref{sec:exps} we compare the confidence bands introduced here with pointwise confidence intervals obtained 
via the asymptotic limit theory as well as with the bootstrap. Of course, pointwise confidence intervals
have a different goal than confidence bands. The pointwise intervals will be be shorter but the
the confidence level will not hold simultaneously across multiple locations.
In contrast, the method we introduce here comes with strong guarantees in terms of finite sample
valid coverage levels across locations.

\section{Constructing a confidence band for a log-concave density}

Given data $X_1,\ldots,X_n$ from a log-concave density $f$ we want to find
functions $\hat{\ell}(x)=\hat{\ell}(x,X_1,\ldots,X_n)$ and $\hat{\mu}(x)=\hat{\mu}(x,X_1,\ldots,X_n)$
such that
$$
\Pr_f \left\{ \hat{\ell}(x) \leq f(x) \leq \hat{\mu}(x)\ \mbox{ for all } x\in \R\right\} \geq 1-\alp
$$
for a given confidence level $1-\alp \in (0,1)$. It is well known that in the case of a general density
$f$ no non-trivial confidence interval $(\hat{\ell}(x),\hat{\mu}(x))$ exists, see e.g. Donoho~(1988).
However, assuming a shape-constraint for $f$ such as log-concavity allows to construct pointwise
and uniform confidence statements as follows:

Let $\CC_n(\alp)$ be a $1-\alp$ confidence set for the distribution function $F$ of $f$, i.e.
\be \label{CS}
\Pr_F \left\{F \in \CC_n(\alp)\right\} \geq 1-\alp.
\ee
Such a nonparametric confidence set always exists, e.g.
the Kolmogorov-Smirnov bands give a confidence set for $F$ (albeit a non-optimal one). Define
\be \label{ell}
\hat{\ell}(x):=\ \inf_{f \mbox{ is log-concave and } F \in \CC_n(\alp)} f(x)
\ee
and define $\hat{\mu}(x)$ analogously with sup in place of inf. If $f$ is log-concave then (\ref{CS})
and (\ref{ell}) imply
$$
\Pr_f \left\{ \hat{\ell}(x) \leq f(x) \leq \hat{\mu}(x)\right\} \geq 1-\alp,
$$
so we obtain a $1-\alp$ confidence interval for $f(x)$ by solving the optimization problem (\ref{ell}).
Moreover, if we solve (\ref{ell}) at multiple locations $x_1,\ldots,x_m$, then we obtain
\be  \label{ells}
\Pr_f \left\{ \hat{\ell}(x_i) \leq f(x_i) \leq \hat{\mu}(x_i)\ \mbox{ for all } i=1,\ldots,m \right\} \geq 1-\alp.
\ee
So the coverage probability is automatically simultaneous across multiple locations and comes with a finite
sample guarantee, since it is inherited from
the confidence set $\CC_n$ in (\ref{CS}). Likewise, the quality of the confidence band, as measured e.g.
by the width $\hat{\mu}(x)-\hat{\ell}(x)$, will also derive from $\CC_n$, which therefore plays a central
role in this approach.
Finally, the log-concavity constraint allows to extend the confidence set (\ref{ells}) 
to a confidence band on the real line, as we will show in Section~\ref{bands}.

Hengartner and Stark~(1995), D\"{u}mbgen~(1998), D\"{u}mbgen~(2003) and Davies and Kovac~(2004) employ the above approach
for inference about a unimodal or a $k$-modal density. Here we introduce a new confidence set $\CC_n(\alp)$ for $F$.
This confidence set is adapted from methodology developed in the abstract Gaussian White Noise model by Walther 
and Perry~(2019) for optimal inference in settings related to the one considered here. Therefore this confidence
set should also prove useful for the works about inference in the unimodal and $k$-modal setting cited above.

The key conceptual problem for solving the optimization problem (\ref{ell}) is that $f$ is infinite dimensional.
We will overcome this by using the log-concavity of $f$ to relax $\CC_n(\alp)$ to a finite dimensional superset,
which makes it possible to compute (\ref{ell}) with fast optimization algorithms.
We will address these tasks in turn in the following subsections.

\subsection{A confidence set for $F$} \label{Cn}

Given $X_1,\ldots,X_n$ i.i.d. from the continuous cdf $F$ we set $s_n:=\lceil \log_2 \log n\rceil$ and
\be  \label{x}
x_i:= X_{\left(1+(i-1)2^{s_n}\right)},\ \ i=1,\ldots,m:=\left\lfloor \frac{n-1}{2^{s_n}} \right\rfloor +1,
\ee
where $X_{(j)}$ denotes the $j$th order statistic.
Our analysis will use only the subset $\{x_i\}$ of the data, i.e. the set containing
every $\log n\,$th order statistic; see Remark~3 for why this is sufficient.

Translating the methodology of the `Bonferroni scan' developed in Walther and Perry~(2019) from the 
Gaussian White Noise model to the density setting suggests employing a confidence set of the form
$$
\CC_n(\alp)\ :=\ \left\{F:\ c_{jkB}\leq F(x_k)-F(x_j) \leq d_{jkB}\mbox{ for all } (j,k)\in 
\II=\bigcup_B \II_B\right\}.
$$
with $c_{jkB}, d_{jkB}, \II$ given below. 
The idea is to choose an index set $\II$ that is rich enough to detect relevant deviations from the empirical
distribution, but which is also sparse enough so that the $|\II|$ constraints can be combined
with a simple weighted Bonferroni adjustment and still result in optimal inference. The second key
ingredient of this construction is to let the weights of the Bonferroni adjustment depend on $j-i$
in a certain way. See Walther and Perry~(2019) for a comparison of the
finite sample and asymptotic performance of this approach with other relevant calibrations, such
as the ones used in the works cited above.

Note that the confidence set $\CC_n(\alp)$ checks the probability content of random
intervals $(x_j,x_k)$, which automatically adapt to the empirical distribution. This makes
it possible to detect relevant deviations from the empirical distribution with a relatively small
number of such intervals, which is key for making the Bonferroni adjustment powerful as well as
 for efficient computation. Moreover, using such random intervals makes the bounds $c_{jkB}, d_{jkB}$
distribution-free since $F(x_k)-F(x_j) \sim \textrm{ Beta}((k-j)2^{s_n}, n+1 -(k-j)2^{s_n})$, see
Ch. 3.1 in Shorack and Wellner~(1986).

The precise specifications of $c_{jkB}, d_{jkB}, \II$ are as follows:

\begin{align*}
\II &:= \bigcup_{B=0}^{B_{max}}\II_B, \ \mbox{ where $B_{max}:=\lfloor \log_2 \frac{n}{8}\rfloor -s_n$}\\
\II_B &:=\Bl\{ (j,k):\ j=1+(i-1)2^B, \ k=1+i2^B \ \mbox{ for } i=1,\ldots,n_B:=\Big\lfloor
\frac{n-1}{2^{B+s_n}} \Big\rfloor \Br\}\\
c_{jkB} &:= c_B:= q\textrm{Beta }\Bl(\frac{\alp}{2(B+2)n_B t_n},2^{B+s_n},n+1-2^{B+s_n}\Br)\\ 
d_{jkB} &:= d_B:= q\textrm{Beta }\Bl(1-\frac{\alp}{2(B+2)n_B t_n},2^{B+s_n},n+1-2^{B+s_n}\Br)
\end{align*}
where $t_n:=\sum_{B=0}^{B_{max}} \frac{1}{B+2}$ and 
$q\textrm{Beta } (\alp, r,s)$ denotes the $\alp$-quantile of the beta distribution with
parameters $r$ and $s$. The term $\frac{1}{B+2}$ is a weighting factor
in the Bonferroni adjustment which results in an advantageous statistical performance, see
Walther and Perry~(2019). It follows from the union bound that $\Pr_F (F \in \CC_n(\alp))
\geq 1-\alp$ whenever $F$ is continuous.

{\bf Remarks: 1.} An alternative way to control the distribution of $F(x_k)-F(x_j)$ is via
a log-likelihood ratio type transformation and Hoeffding's inequality, see Rivera and Walther~(2013)
and Li, Munk, Sieling and Walther~(2020). This results in a loss of power due to the slack in 
Hoeffding's inequality and the slack from inverting the log-likelihood ratio transformation with
an inequality. Simulations show that the above approach using an exact beta distribution is less
conservative despite the use of Bonferroni's inequality to combine the statistics across $\II$.

{\bf 2.} The inference is based on the statistic $F((x_j,x_k))$, i.e. the unknown $F$ evaluated on
the random interval $(x_j,x_k)$, rather than on the more commonly used statistic $F_n(I)$, which
evaluates the empirical measure on deterministic intervals $I$. The latter statistic follows a 
binomial distribution whose discreteness makes it difficult to combine these statistics
across $\II$ using Bonferroni's inequality without incurring substantial conservatism
and hence loss of power.
This is another important reason for using random intervals $(x_j,x_k)$ besides the adaptivity
property mentioned above. Moreover, a deterministic system of intervals would have to be anchored
around the range of the data and this dependence on the data is difficult to account for and is therefore
typically glossed over in the inference.

{\bf 3.} The definition of $x_i$ in (\ref{x}) means that we do  not consider intervals $(X_{(j)},X_{(k)})$
with $k-j<\log n$. Thus, as opposed to the regression setting in Walther and Perry~(2019) we omit the first
block\footnote{We also shift the index $B$ to let it start at 0 rather than at 2. This results in
a simpler notation but does not change the methodology.} of intervals. This derives from the folklore knowledge
in density estimation that at least $\log n$ observations are required in order to obtain consistent inference
simultaneously across such intervals. Indeed, this choice is sufficient to yield the 
asymptotic optimality result in  Theorem~\ref{thm1}.

We further simplify the construction in Walther and Perry~(2019) by restricting ourselves to a dyadic spacing of the indices
$k-j$ since we already obtain quite satisfactory results with this set of intervals. 

\subsection{Bounds for $\int_a^b f$ when $f$ is log-concave}

The confidence set $\CC_n(\alp)$ describes a set of plausible distributions in terms of $\int_{x_j}^{x_k} f(t)\, dt$
for infinite dimensional $f$. In the special case when $f$ is log-concave it is possible to construct
a finite dimensional superset of $\CC_n(\alp)$ by deriving bounds for this integral
in terms of functions of a finite number of variables:

\begin{Lemma} \label{integralbounds}
Let $f$ be a univariate log-concave function. For given $x_1 <\ldots <x_m$ write $\ell_i:=\log f(x_i)$,
$i=1,\ldots,m$. Then there exist real numbers $g_2,\ldots,g_{m-1}$ such that
$$
\ell_j \leq \ell_i +g_i \,(x_j-x_i)\ \ \mbox{ for all } i \in \{2,\ldots,m-1\},\ j\in \{i-1,i+1\}
$$
and 
$$
(x_{i+1}-x_i) \exp(\ell_{i}) \,E(\ell_{i+1}- \ell_i)\ \leq \
\int_{x_i}^{x_{i+1}} f(t)\,dt \ \leq \ 
\begin{cases}
\exp (\ell_i) (x_{i+1}-x_i)\, E \bigl(g_i(x_{i+1}-x_i)\bigr),  &i \in \{2,\ldots,m-1\}\\
\exp (\ell_{i+1}) (x_{i+1}-x_i) \, E \bigl(g_{i+1}(x_i-x_{i+1})\bigr), &i \in \{1,\ldots,m-2\}
\end{cases}
$$
where
$$
E(s)\ :=\ \int_0^1 \exp(st)\,dt\ =\begin{cases}
\frac{\exp(s) -1}{s} & \mbox{ if } s \neq 0\\
1 & \mbox{ if } s = 0
\end{cases}
$$
is a strictly convex and infinitely often differentiable function.
\end{Lemma}
Importantly, the bounds given in the lemma are convex and smooth functions of the $g_i$ and $\ell_i$, despite
the fact that these variables appear in the denominator in the formula for $E$. 
This makes it possible to bring fast optimization
routines to bear on the problem (\ref{ell}).

\subsection{Computing pointwise confidence intervals }

We are now in position to define a superset of $\CC_n(\alp)$ by relaxing the inequalities
$c_{B}\leq \int_{x_j}^{x_k} f(t)\, dt \leq d_{B}$ in the definition of $\CC_n(\alp)$.
To this end define for $i=1,\ldots,m-1$ the functions
\begin{align*}
L_i(\x,\bsell) &:= (x_{i+1}-x_i) \exp(\ell_i)\,E(\ell_{i+1}- \ell_i) \\
U_i(\x,\bsell,\g) &:= \begin{cases}
\exp (\ell_i) (x_{i+1}-x_i) \,E\bigl(g_i(x_{i+1}-x_i)\bigr),  &i =m-1\\
\exp (\ell_{i+1}) (x_{i+1}-x_i) \,E\bigl(g_{i+1}(x_i-x_{i+1})\bigr), &i \in \{1,\ldots,m-2\}
\end{cases}  \\
V_i(\x,\bsell,\g) &:= \begin{cases}
\exp (\ell_i) (x_{i+1}-x_i) \,E \bigl(g_i(x_{i+1}-x_i)\bigr),  &i \in \{2,\ldots,m-1\}\\
\exp (\ell_{i+1}) (x_{i+1}-x_i)\,E \bigl(g_{i+1}(x_i-x_{i+1})\bigr), &i =1
\end{cases} 
\end{align*}
where $\x=(x_1,\ldots,x_m)$, $\bsell=(\ell_1,\ldots,\ell_m)$, $\g=(g_1,\ldots,g_{m-1})$ and $E(\cdot)$
is defined in Lemma~\ref{integralbounds}.

Given the subset $x_1 < \ldots < x_m$ of the order statistics defined in (\ref{x}), we define
$\tilde{\CC}_n(\alp)$ to be the set of densities $f$ for which there exist real $g_1,\ldots,g_{m-1}$
such that $\ell_i:=\log f(x_i)$,
$i=1,\ldots,m$, satisfy (\ref{CONC})-(\ref{UP}):
\begin{align}
\ell_j & \leq \ell_i +g_i (x_j-x_i)\ \ \mbox{ for all } i \in \{2,\ldots,m-1\},\ j\in \{i-1,i+1\} 
\label{CONC} \\
\mbox{For }B=0,\ldots,B_{max}&: \nn \\
c_{B}  \ &\leq \ \sum_{i=j}^{k-1} U_i(\x,\bsell,\g)\ \ \mbox{ for all } (j,k)\in
\II_B \label{DOWN1}\\
c_{B}  \ &\leq \ \sum_{i=j}^{k-1} V_i(\x,\bsell,\g)\ \ \mbox{ for all } (j,k)\in
\II_B \label{DOWN2}\\
\sum_{i=j}^{k-1} L_i(\x,\bsell)\ & \leq \ d_{B} \ \ \mbox{ for all } (j,k)\in
\II_B. \label{UP}
\end{align}

Now we can implement a computable version of the confidence bound (\ref{ell})
by optimizing over $\tilde{\CC}_n(\alp)$ rather than over $\CC_n(\alp)$. Note that if $f$
is log-concave then it follows from Lemma~\ref{integralbounds} that $f \in \CC_n(\alp)$
implies $f \in \tilde{\CC}_n(\alp)$. This proves the following key result:

\begin{Proposition}  \label{valid}
If $f$ is log-concave then $\Pr_f \{f \in \tilde{\CC}_n(\alp)\} \geq 1-\alp$. 
Consequently, if we define pointwise lower and upper confidence bounds at the $x_i$, $i=1,\ldots,m$, 
via the optimization problems
\begin{align}
\hat{\ell}(x_i):=\ &\min \ \ell (x_i) \label{Opt}\\
& \mbox{subject to } f \in \tilde{\CC}_n(\alp), \mbox{ i.e. subject to (\ref{CONC})-(\ref{UP})} \nn \\
\hat{\mu}(x_i):=\ &\max \ \ell (x_i) \nn \\
& \mbox{subject to } f \in \tilde{\CC}_n(\alp), \mbox{ i.e. subject to (\ref{CONC})-(\ref{UP})}, \nn
\end{align}
then the following simultaneous confidence statement holds whenever $f$ is log-concave:
$$
\Pr_f \left\{ \exp (\hat{\ell}(x_i)) \leq f(x_i) \leq \exp (\hat{\mu}(x_i))
\ \mbox{ for all } i=1,\ldots,m \right\} \geq 1-\alp.
$$
\end{Proposition}

It is an important feature of these confidence bounds that they come with a finite sample guaranteed
confidence level $1-\alp$. On the other hand, it is desirable that the construction is not overly
conservative (i.e. has coverage not much larger than $1-\alp$) as otherwise it would result in unnecessarily 
wide confidence bands. This is the motivation for deriving a  statistically optimal confidence set
in Section~\ref{Cn} and for deriving bounds in Lemma~\ref{integralbounds} that are sufficiently
tight.
Indeed, it will be shown in Section~\ref{asym}  that the above construction results in statistically
optimal confidence bands.

\subsection{Constructing confidence bands} \label{bands}

The simultaneous pointwise confidence bounds $\left(\hat{\ell}(x_i), \hat{\mu}(x_i)\right), i=1,\ldots,m$,
from the optimization problem (\ref{Opt}) imply a confidence band on the real line due to
the concavity of $\log f$. 
In more detail, we can extend the definition of $\hat{\ell}$ to the real line
simply by interpolating between the $\hat{\ell}(x_i)$:
\be \label{lower}
\hat{\ell}(x) := \begin{cases}
\hat{\ell}(x_i)+(x-x_i)\frac{\hat{\ell}(x_{i+1}) - \hat{\ell}(x_i)}{x_{i+1}-x_i} & \mbox{if }
x\in [x_i,x_{i+1}),\ i\in \{1,\ldots,m-1\} \\
-\infty & \mbox{otherwise.}
\end{cases}
\ee
Then $\log f(x_i) \geq \hat{\ell}(x_i)$ for $i=1,\ldots,m$ implies $\log f(x) \geq \hat{\ell}(x)$ for $x \in \R$
since $\log f$ is concave and $\hat{\ell}$ is piecewise linear. (In fact, it follows from (\ref{Opt})
that $\hat{\ell}$ is also concave.)

In order to construct an upper confidence bound note that concavity of $\log f$ 
together with $\hat{\ell}(x_i) \leq \log f(x_i) \leq \hat{\mu}(x_i)$ for all $i \in \{1,\ldots,m\}$ implies for $x>x_k$
with $k\in \{2,\ldots,m\}$:
\begin{align*}
\frac{\log f(x) -\hat{\mu}(x_k)}{x-x_k} & \leq \frac{\log f(x) -\log f(x_k)}{x-x_k}
\leq \min_{j\in \{1,\ldots,k-1\}} \frac{\log f(x_k) -\log f(x_j)}{x_k-x_j}\\
& \leq \min_{j\in \{1,\ldots,k-1\}} \frac{\hat{\mu}(x_k) -\hat{\ell}(x_j)}{x_k-x_j} \ =:L_k
\end{align*}
and likewise for $x<x_k$ with $k\in \{1,\ldots,m-1\}$:
$$
\frac{\hat{\mu}(x_k)-\log f(x)}{x_k-x}\, \geq \, 
\max_{j\in \{k+1,\ldots,m\}} \frac{\hat{\ell}(x_j) -\hat{\mu}(x_k)}{x_j-x_k}\ =:R_k
$$
Hence $\log f$ is bounded above by
$$
\hat{\mu}(x) :=
\begin{cases}
\hat{\mu}(x_{i+1})+R_{i+1} (x-x_{i+1}) & \mbox{if }x\in (x_i,x_{i+1}), \ i \in \{0,1\}, \ \mbox{ where } x_0:=-\infty\\
M_i(x)  & \mbox{if }x\in [x_i,x_{i+1}), \ i \in \{2,\dots, m-2\} \\
\hat{\mu}(x_i)+L_i(x-x_i) & \mbox{if }x\in [x_i,x_{i+1}), \ i \in \{m-1,m\}, \ \mbox{ where } x_{m+1}:=\infty
\end{cases}
$$
with
\begin{align*}
M_i(x) & := \min \Bl( \hat{\mu}(x_i)+L_i(x-x_i), \,\hat{\mu}(x_{i+1})+R_{i+1} (x-x_{i+1})\Br) \\
& = \Bl( \hat{\mu}(x_i)+L_i(x-x_i)\Br) 1 (x \in [x_i,\bar{x}_i))\, +\,
\Bl(\hat{\mu}(x_{i+1})+R_{i+1} (x-x_{i+1})\Br) 1 (x \in [\bar{x}_i,x_{i+1}))
\end{align*}
where $\bar{x}_i:= \frac{\hat{\mu}(x_{i+1})-\hat{\mu}(x_i)+L_i x_i -R_{i+1} x_{i+1}}{L_i-R_{i+1}}$.

Thus we proved:

\begin{Proposition}  \label{valid2}
If $f$ is log-concave then
$$
\Pr_f \left\{ \exp (\hat{\ell}(x)) \leq f(x) \leq \exp (\hat{\mu}(x))
\ \mbox{ for all } x \in \R \right\} \geq 1-\alp.
$$
\end{Proposition}

The upper bound $\hat{\mu}(x)$ need not be concave but it is minimal in the sense that it can be shown that
for every real $x$ there exist a concave function $g$ with $\hat{\ell}(x_i) \leq g(x_i) \leq \hat{\mu}(x_i)$
for all $i \in \{1,\ldots,m\}$ and $g(x)=\hat{\mu}(x)$.

As an alternative to $\exp(\hat{\mu}(x))$ we tried a simple interpolation between the points $(x_i,\exp(\hat{\mu}(x_i)))$.
This interpolation will result in a smoother bound than $\exp(\hat{\mu}(x))$, but the coverage guarantee of
Proposition~\ref{valid2} does not apply any longer. However, the difference between $\exp(\hat{\mu}(x))$
and the interpolation will vanish as the sample size increases (or by increasing the number $m$ of design points
$x_i$ in (\ref{x}) for a given sample size $n$), and the simulations in Section~\ref{sec:exps}
show that the empirical coverage exceeds the nominal level in all cases considered. Therefore we
also recommend this interpolation as a simple and smoother alternative to $\exp(\hat{\mu}(x))$. 

Finally, we point out that the computational effort can be lightened simply by solving the
optimization problem (\ref{Opt}) for a subset of $\{x_i,\,1 \leq i \leq m\}$ and then constructing
$\hat{\ell}$ and $\hat{\mu}$ as described above based on that smaller subset of $x_i$. Such a confidence band
will still satisfy Proposition~\ref{valid2}, but it will be somewhat wider at locations between those design
points $x_i$ as it is based on fewer pointwise confidence bounds.
Hence there is a trade-off between the width of the band and the computational effort required. 
While a larger subset of the $x_i$ will result in a somewhat reduced width of the band between the $x_i$,
there are diminishing returns as the width at the $x_i$ will not change.
It follows from Theorem~\ref{thm1} below that solving the optimization problem (\ref{Opt}) for 
$\{x_i,\,1 \leq i \leq m\}$ with 
$m$ given in (\ref{x}) is sufficient to produce statistically optimal confidence bands in a representative
setting.

\section{Solving the optimization problem}

Next we describe a method for computing the pointwise confidence intervals $(\hat \ell(x_i), \hat \mu(x_i)), \; 
i=1,\ldots,m$, from observations $X_i, \; i=1,\ldots,n$, by efficiently solving the optimization problems \eqref{Opt}.  
Constructing the confidence band is then straightforward with
the post-processing steps given in Section \ref{bands}.

Inspecting the optimization problems \eqref{Opt}, we see that these problems possess some interesting structure: 
The criterion functions are linear, and the constraints \eqref{CONC} and \eqref{UP} are convex.  However, the 
constraints \eqref{DOWN1} and \eqref{DOWN2} are non-convex. Finding the global minimizer of a non-convex optimization 
problem (even a well-structured one) can be challenging; instead, we focus on a method for finding critical points 
of the problems \eqref{Opt}.  Taking a closer look at the non-convex constraints \eqref{DOWN1} and \eqref{DOWN2}, we 
make the simple observation that they may be expressed as the difference of two convex functions (namely, a constant 
function minus a convex function).  This property puts the problems \eqref{Opt} into the special class of non-convex 
problems commonly referred to as \textit{difference of convex programs} (Hartman, 1959; Tao, 1986; Horst and 
Thoai, 1999; Horst et al., 2000), for which a critical point can be efficiently found.  The class of difference of 
convex programs is quite broad, encompassing many problems encountered in practice, with a good amount of research 
into this area continuing on today. Important references include Hartman (1959), Tao (1986), Horst and Thoai (1999), 
Horst et al.~(2000), Yuille and Rangarajan (2003), Smola et al.~(2005), and Lipp and Boyd (2016).

% \clearpage
\begin{algorithm}[h!]
    \caption{Penalty convex-concave procedure for computing pointwise confidence intervals for a log-concave density}
    \label{alg:ccp}
    \begin{algorithmic}
        \STATE {\bf Input:} Subset of observations $x_t, \; t=1,\ldots,m$; collection of interval endpoints 
$\II=\bigcup_B\II_B$; initial points $\bsell^{(0)}$, $\g^{(0)}$; initial penalty strength $\tau_0 > 0$; 
penalty growth factor $\kappa > 1$; maximum penalty strength $\tau_{\max} > \tau_0$; maximum number of iterations $K_{\max}$
        
        \STATE {\bf For $t=1,\ldots,m$, $K=0,\ldots,K_{\max}$:} \\
        
        \STATE \quad Convexify the constraints \eqref{DOWN1} and \eqref{DOWN2}, by forming the linearizations, for $i=1,\ldots,m-1$,
        \begin{align*}
        \hat U_i(\x, \bsell, \g; \bsell^{(K)}, \g^{(K)}) & = U_i(\x, \bsell^{(K)}, \g^{(K)}) + \big\langle \nabla U_i(\x, \bsell^{(K)}, \g^{(K)}) ,\, (\bsell, \g) - (\bsell^{(K)}, \g^{(K)}) \big\rangle \\ %\label{eq:DOWN1_lin} \\
        \hat V_i(\x, \bsell, \g; \bsell^{(K)}, \g^{(K)}) & = V_i(\x, \bsell^{(K)}, \g^{(K)}) + \big\langle \nabla V_i(\x, \bsell^{(K)}, \g^{(K)}) ,\, (\bsell, \g) - (\bsell^{(K)}, \g^{(K)}) \big\rangle. %\label{eq:DOWN2_lin}
        \end{align*}

        \STATE \quad Solve the pair of convexified problems
        \begin{align}
            \ell_t^{(K+1)} = & \min \ \ell_t + \tau_{K} \cdot \sum_{B=0}^{B_{max}} \sum_{(j,k) \in \II_B} s_{B,j,k}
                        \label{eq:ccp:min} \\
            & \mbox{subject to \eqref{CONC}, \eqref{UP}}, \textrm{ and} \notag \\
            & \hspace{0.65in} \mbox{For }B=0,\ldots,B_{max}: \nn \notag \\
            & \hspace{1.75in} c_{B} - \sum_{i=j}^{k-1} \hat U_i(\x, \bsell, \g; \bsell^{(K)}, \g^{(K)}) \leq s_{B,j,k} \; \mbox{ for all } (j,k) \in \II_B \notag \\
            & \hspace{1.75in} c_{B} - \sum_{i=j}^{k-1} \hat V_i(\x, \bsell, \g; \bsell^{(K)}, \g^{(K)}) \leq s_{B,j,k} \; \mbox{ for all } (j,k) \in \II_B \notag \\
            & \hspace{1.75in} s_{B,j,k} \geq 0 \; \mbox{ for all } (j,k) \in \II_B \notag
        \end{align}
        \begin{align}
            \mu_t^{(K+1)} = & \max \ \ell_t - \tau_{K} \cdot \sum_{B=0}^{B_{max}} \sum_{(j,k) \in \II_B} s_{B,j,k}
                 \label{eq:ccp:max} \\
            & \mbox{subject to \eqref{CONC}, \eqref{UP}}, \textrm{ and} \notag\\
            & \hspace{0.65in} \mbox{For }B=0,\ldots,B_{max}: \nn \notag \\
            & \hspace{1.75in} c_{B} - \sum_{i=j}^{k-1} \hat U_i(\x, \bsell, \g; \bsell^{(K)}, \g^{(K)}) \leq s_{B,j,k} \; \mbox{ for all } (j,k) \in \II_B \notag\\
            & \hspace{1.75in} c_{B} - \sum_{i=j}^{k-1} \hat V_i(\x, \bsell, \g; \bsell^{(K)}, \g^{(K)}) \leq s_{B,j,k} \; \mbox{ for all } (j,k) \in \II_B \notag \\
            & \hspace{1.75in} s_{B,j,k} \geq 0 \; \mbox{ for all } (j,k) \in \II_B. \notag
        \end{align}        

        \STATE \quad Update the penalty strength, by setting $\tau_{K+1} = \min \{ \kappa \cdot \tau_{K},\, \tau_{\max} \}$.

        \STATE {\bf Output:} Pointwise confidence intervals $(\ell_t^{(K)} \mu_t^{(K)}), \; t=1,\ldots,m.$
    \end{algorithmic}
\end{algorithm}
% \clearpage

A natural approach to finding a critical point of a difference of convex program is to linearize the non-convex 
constraints, then solve the convexified problem using any suitable off-the-shelf solver, and repeating these steps as 
necessary.  This strategy underlies the well-known convex-concave procedure in Yuille and Rangarajan~(2003), a popular 
heuristic for difference of convex programs. In more detail, the convex-concave iteration as applied to the problems 
\eqref{Opt} works as follows:  Given feasible initial points, we first replace the (non-convex) constraints \eqref{DOWN1}
and \eqref{DOWN2} by their first-order Taylor approximations centered around the inital points.  Formally, letting 
$\bsell^{(K)}$ and $\g^{(K)}$ denote the log-densities and subgradients on iteration $K$, respectively, we form 
% the linearizations
\begin{align}
\hat U_i(\x, \bsell, \g; \bsell^{(K)}, \g^{(K)}) & := U_i(\x, \bsell^{(K)}, \g^{(K)}) + \big\langle \nabla U_i(\x, \bsell^{(K)}, \g^{(K)}) ,\, (\bsell, \g) - (\bsell^{(K)}, \g^{(K)}) \big\rangle \label{eq:DOWN1_lin} \\
\hat V_i(\x, \bsell, \g; \bsell^{(K)}, \g^{(K)}) & := V_i(\x, \bsell^{(K)}, \g^{(K)}) + \big\langle \nabla V_i(\x, \bsell^{(K)}, \g^{(K)}) ,\, (\bsell, \g) - (\bsell^{(K)}, \g^{(K)}) \big\rangle, \label{eq:DOWN2_lin}
\end{align}
for $i=1,\ldots,m-1$.  We then solve the convexified problems (using the constraints \eqref{eq:DOWN1_lin}, 
\eqref{eq:DOWN2_lin} instead of \eqref{DOWN1}, \eqref{DOWN2}) with any off-the-shelf solver.  Then we re-compute the 
approximations using the obtained solutions to the convexified problems and repeat these steps until an appropriate 
stopping criterion has been satisfied (e.g., some pre-determined maximum number of iterations has been reached, the 
change in criterion values are smaller than some specified tolerance, the sum of the slack variables is less than some tolerance, and/or we have that $\tau_K \geq \tau_{\max}$).  From this description, it may be apparent to 
the reader that the convex-concave procedure is actually a generalization of the majorization-minimization class of 
algorithms (which includes the well-known expectation-maximization algorithm as a special case).

We give a complete description of the convex-concave procedure as applied to the optimization problems \eqref{Opt} in 
Algorithm \ref{alg:ccp} appearing above, along with one important modification that we explain now.  In practice
 it is not easy to obtain feasible initial points for the problems \eqref{Opt}.  Therefore, the penalty convex-concave 
procedure, a modification to the basic convex-concave procedure that was introduced by Le Thi and Dinh (2014), 
Dinh and Le Thi (2014), and Lipp and Boyd (2016), works around this issue by allowing for an (arbitrary) infeasible 
initial point and then gradually driving the iterates into feasibility by adding a penalty for constraint violations 
into the criterion that grows with the number of iterations (explaining the word ``penalty'' in the name of 
the procedure), through the use of slack variables.

Standard convergence theory for the (penalty) convex-concave procedure (see, e.g., Section 3.1 in Lipp and Boyd (2016) 
as well as Theorem~10 in Sriperumbudur and Lanckriet (2009) and Proposition~1 in Khamaru and Wainwright (2018) tells 
us that the criterion values \eqref{eq:ccp:min} and \eqref{eq:ccp:max} generated by Algorithm \ref{alg:ccp} converge.  
Moreover, under regularity conditions (see Section 3.1 in Lipp and Boyd (2016)), the iterates generated by Algorithm 
\ref{alg:ccp} converge to critical points of the problems \eqref{Opt}.  At convergence, the pointwise confidence
intervals generated by Algorithm~\ref{alg:ccp} can be turned in confidence bands as described in Section~\ref{bands}.

Finally, we mention that although not necessary, additionally linearizing the (convex) constraints \eqref{UP} around 
the previous iterate, i.e., forming
\[
L_i(\x, \bsell^{(K)}) + \big\langle \nabla L_i(\x, \bsell^{(K)}), \, \bsell - \bsell^{(K)} \big\rangle, \; i=j,\ldots,k-1, \; (j,k) \in \II_B, \; B=0,\ldots,B_{max},
\]
on iteration $K$, can help circumvent numerical issues.  Furthermore, as all the constraints in the problems 
\eqref{eq:ccp:min} and  \eqref{eq:ccp:max} are now evidently affine functions, this relaxation has the added benefit of 
turning the problems \eqref{eq:ccp:min} and \eqref{eq:ccp:max} into linear programs, for which there (of course) exist 
heavily optimized solvers.

\section{Large sample statistical performance}  \label{asym}

The large sample performance of the log-concave MLE has been studied intensively, see e.g. D\"{u}mbgen
and Rufibach~(2009), Kim and Samworth~(2016) and Doss and Wellner~(2016). The main message is that
the MLE attains the optimal minimax rate of convergence of $O\left(n^{-2/5}\right)$ with respect to various
global loss functions. Recently, Kim, Guntuboyina and Samworth~(2018) have shown that the MLE
can achieve a faster rate of convergence when the log density is $k$-affine, i.e. when $\log f$ consists of $k$
linear pieces . They show that the MLE is able to adapt to this simpler model, where it will converge
with nearly the parametric rate, namely with $O\left(n^{-1/2} \log^{5/8}n\right)$. Here we show that the
construction of our confidence band via the particular confidence set $\CC_n(\alp)$ will also result
in a nearly parametric rate of convergence for the width of the confidence band in that case. To this end, 
we first consider the case where some part of $f$ is log-linear:

\begin{Theorem} \label{thm1}
Let $f$ be a log-concave density and suppose $\log f$ is linear on some interval $J$. (So $J$ may be
a proper subset of the support of $f$). 
Then on every closed interal $I \subset int\,J$:

$$
\Pr_f \left\{ \sup_{x \in I} \Bl|\hat{u}(x) -f(x)\Br| 
\ \leq \ C \sqrt{\frac{ \log \log n}{n}} \right\} \ra 1
$$

for some constant $C$, and the same statement holds for $\hat{\ell}$ in place of $\hat{u}$.
 In particular, the width of the confidence band satisfies 
$\max_{x \in I} \bigl(\hat{u}(x)-
\hat{\ell}(x)\bigr) \leq 2C \sqrt{\frac{\log \log n}{n}}$ with probability converging to 1.
\end{Theorem}

If there are $k$ such intervals, then the theorem holds for the maximum width over the $k$ intervals.
This includes $k$-affine log densities as a special case.

We conjecture that the width of the confidence band will likewise achieve the optimal minimax rate
if $\log f$ is smooth rather than linear.

\section{Some examples}
\label{sec:exps}

Finally, we present some numerical examples of our methodology, highlighting the empirical coverage and widths of our 
confidence bands as well as the computational cost of computing the bands, for a number of different distributions.

To calculate coverage, we first simulated $n \in \{100,1000\}$ observations from a (i) standard normal distribution, 
(ii) uniform distribution on $[-10,10]$, (iii) chi-squared distribution with three degrees of freedom, and 
(iv) exponential distribution with parameter 1.  Then, we computed our confidence bands from the data by 
running the penalty convex-concave procedure described in Algorithm \ref{alg:ccp} and then computing
$\exp(\hat{\ell}(x))$ and $\exp(\hat{\mu}(x))$ as discussed in Section~\ref{bands}, where $\exp(\hat{\mu}(x))$
was computed by linearly interpolating between the $(x_i,\exp(\hat{\mu}(x_i)))$.
We repeated these two steps (simulating data and computing bands) 1000 times. We calculated the empirical 
coverage for each density $f$ by evaluating the band at 10000 points $\{t_j\}$, evenly spaced across the range of the data,
to check whether $\exp(\hat{\ell}(t_j)) \leq f(t_j) \leq \exp(\hat{\mu}(t_j))$ for all $j$, and then computed the empirical
frequency of this event across the 1000 repetitions.
In order to calculate the widths of the bands, we averaged the widths at the sample quartiles over all of the 
repetitions. We calculated the computational effort by averaging the runtimes, obtained by running 
Algorithm~\ref{alg:ccp} on a workstation with four Intel E5-4620 2.20GHz processors and 15 GB of RAM, over all  
the repetitions.  To speed up the computation for $n\geq 1000$, 
we ran Algorithm \ref{alg:ccp} on a subset of 30\% of the points 
$x_i, \; i=1,\ldots,m$, as discussed at the end of Section~\ref{bands}; the coverages and widths were virtually 
indistinguishable from those obtained using the full set of points $x_i, \; i=1,\ldots,m$.

Algorithm \ref{alg:ccp} requires a few tuning parameters, which are important for assuring quick convergence.  
In general, we found that the initial penalty strength $\tau_0$, the maximum penalty strength $\tau_{\max}$, and 
the penalty growth factor $\kappa$ had the greatest impact on convergence.  In our experience, setting $\tau_0$ 
to a small value and $\tau_{\max}$ to a large value worked well; we used $\tau_0 \in \{10^{-5}, 10^{-4}, 10^{-3}\}$ 
and $\tau_{\max} \in \{10^{3}, 10^{4}, 10^{5}\}$, with the most suitable values depending on the characteristics of 
the problem.  We experimented with various penalty growth factors $\kappa \in \{1,2,\ldots,10\}$, finding that the 
best value of $\kappa$ again varied with the problem setup.  We always set the maximum number of iterations 
$K_{\max} = 50$ as our method usually converged after around 20-30 iterations across all problem settings.  
We initialized the points $\bsell^{(0)}, \g^{(0)}$ randomly.  
Finally, we set the miscoverage level $\alpha = 0.1$ but we also report results for $\alpha = 0.05$.

\begin{figure*}[h!]
\vskip 0.2in
% \begin{center}
% \centerline{
\centering
\includegraphics[width=0.49\textwidth]{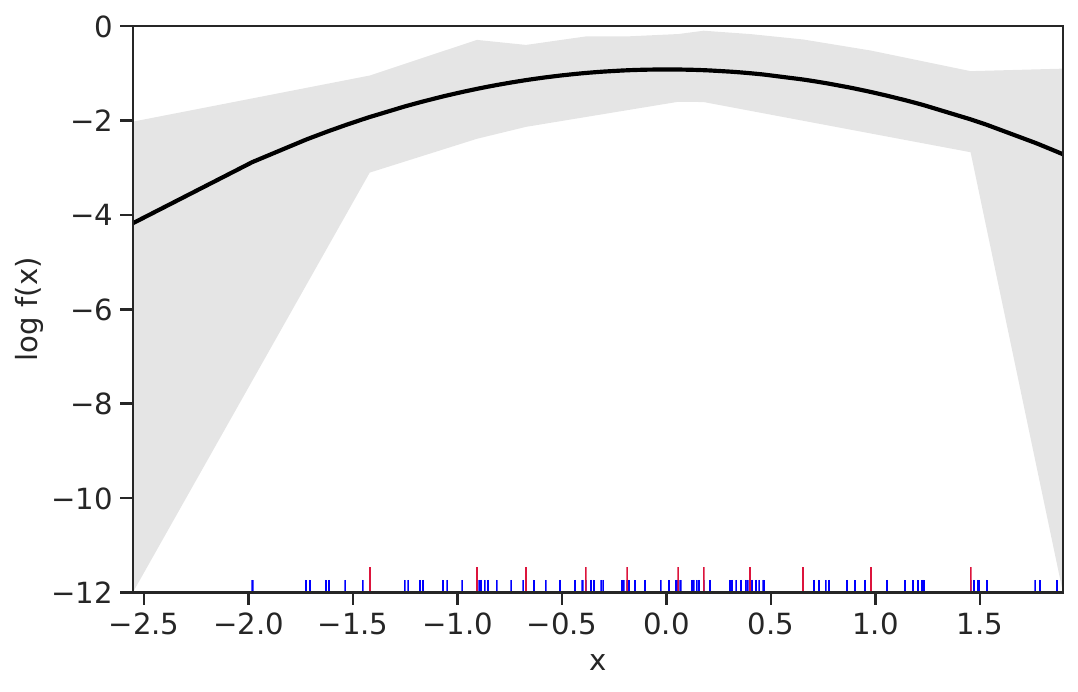} \hfill
\includegraphics[width=0.49\textwidth]{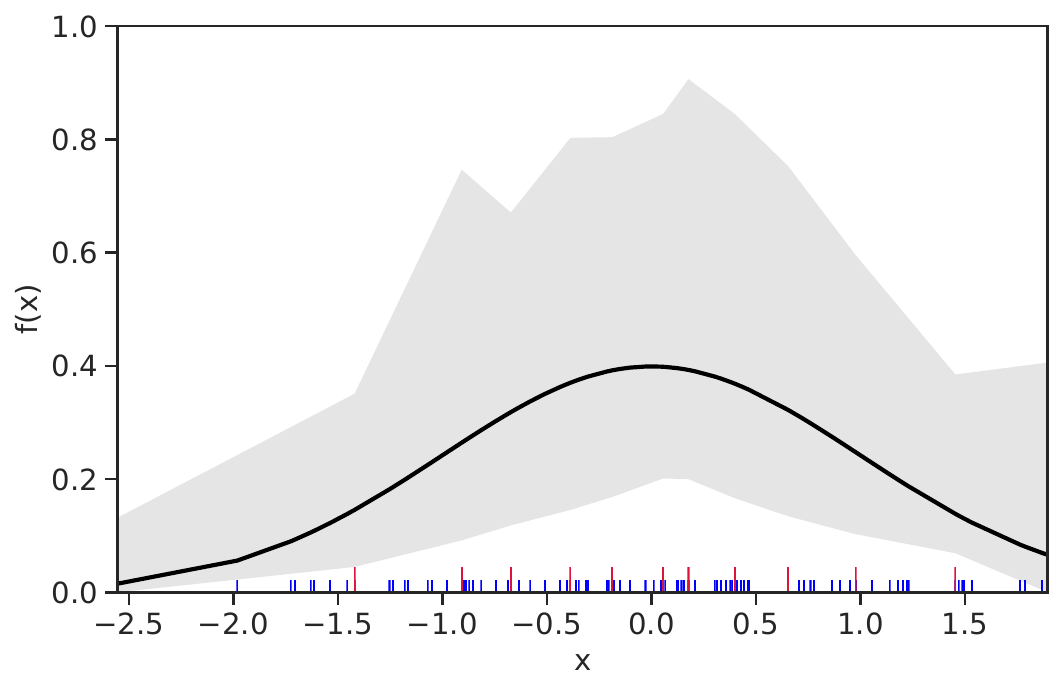} \\
\includegraphics[width=0.49\textwidth]{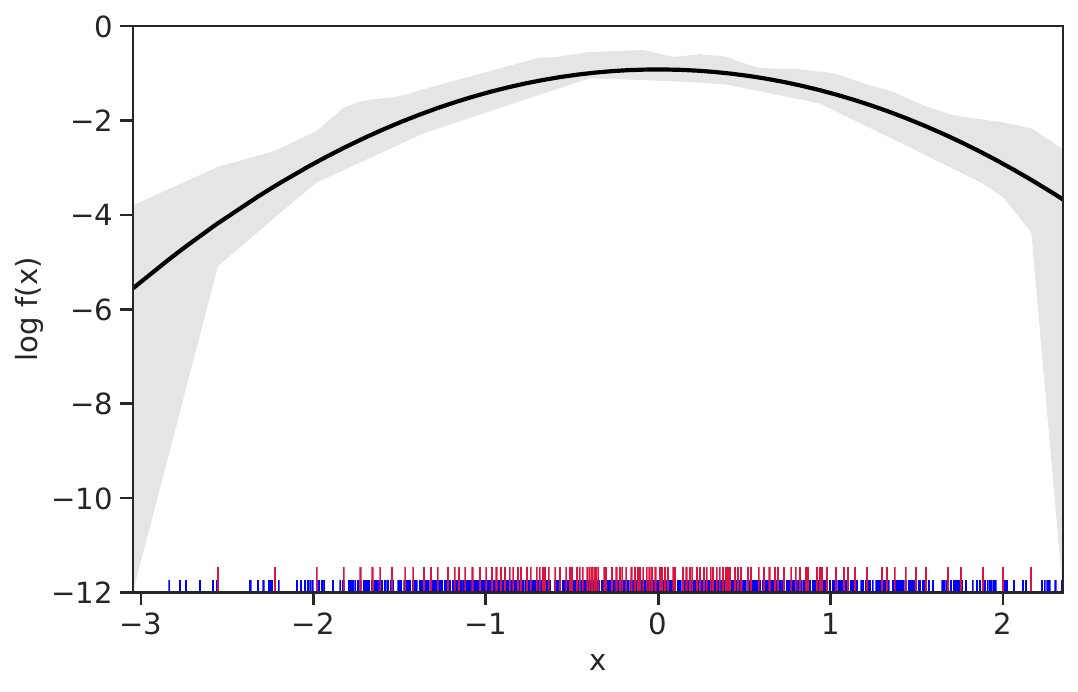} \hfill
\includegraphics[width=0.49\textwidth]{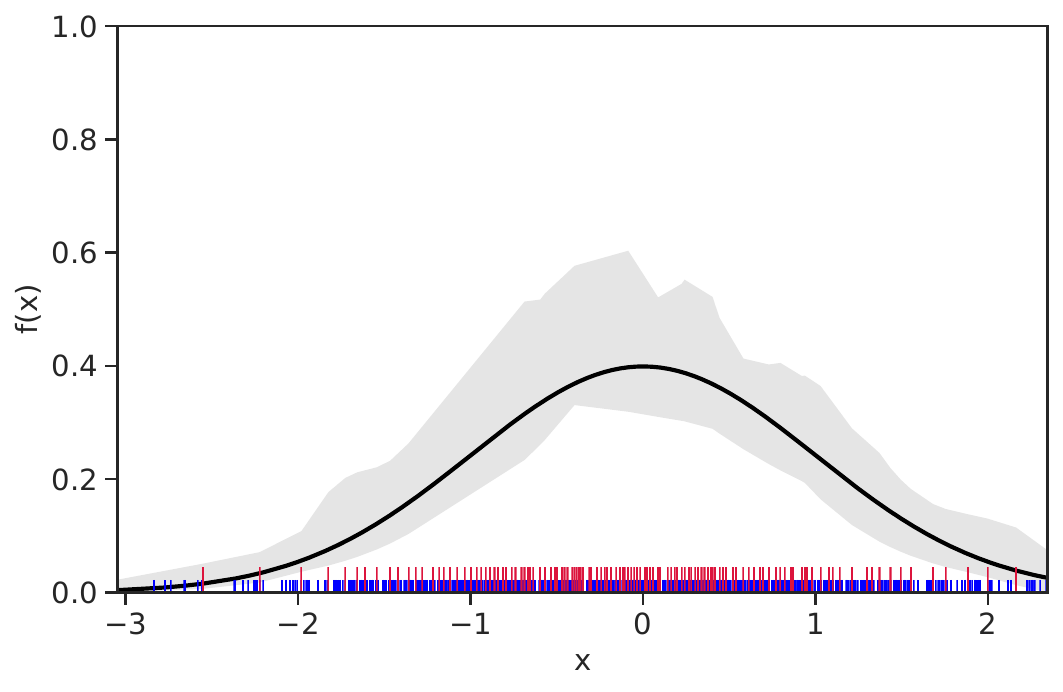} \\
\includegraphics[width=0.49\textwidth]{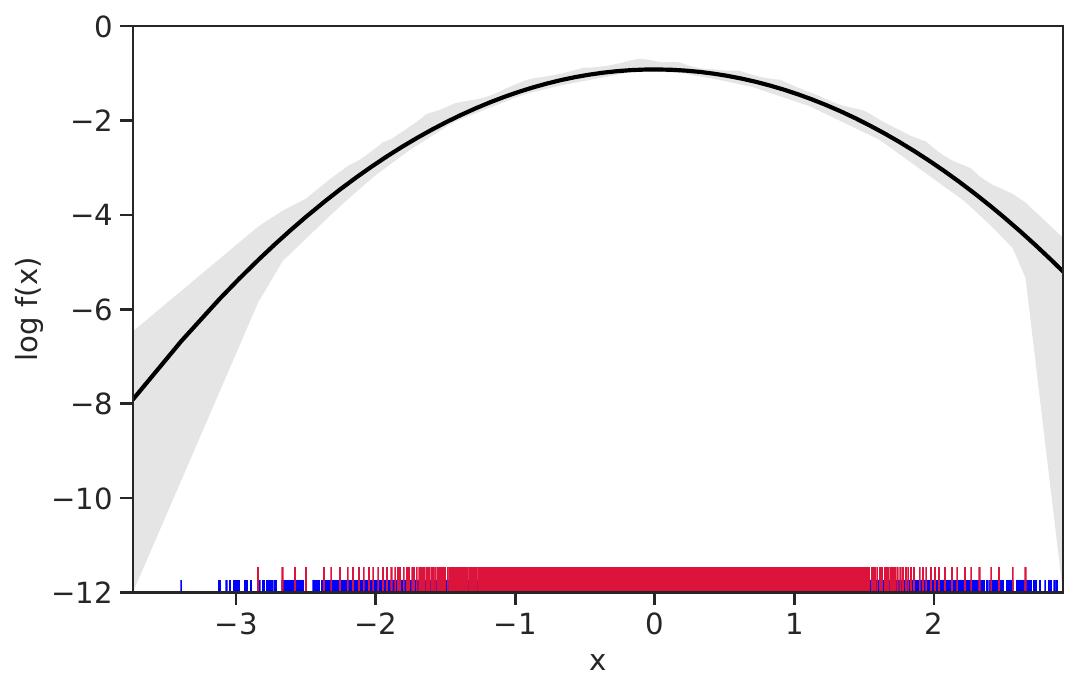} \hfill
\includegraphics[width=0.49\textwidth]{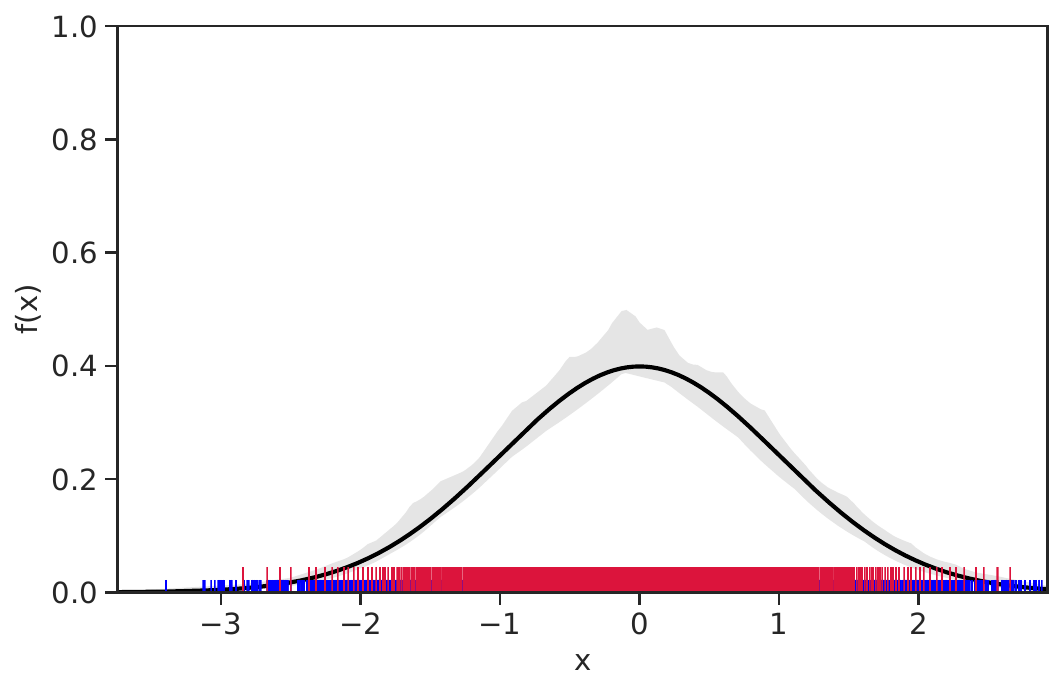}
% }
\vskip -0.2in
\caption{\textit{Confidence bands (shaded in gray) generated by Algorithm \ref{alg:ccp} and linear interpolation  
as described in Section~\ref{bands}, for a Gaussian density.  The left column shows the bands for the log density, 
while the right column shows the bands for the density. The solid black line marks the underlying (log)density.  
Top, middle and bottom rows show results for sample sizes $n= 100$, $1000$, and $10000$, respectively.  
At the bottom of each plot, the observations
$X_i$ are indicated in blue (short lines), while the points $x_i, \; i=1,\ldots,m$ are marked in red (long lines).}}
\label{fig:bands:gaus}
% \end{center}
\vskip -0.2in
\end{figure*}

\begin{figure*}[h!]
\vskip 0.2in
% \begin{center}
% \centerline{
\centering
\includegraphics[width=0.49\textwidth]{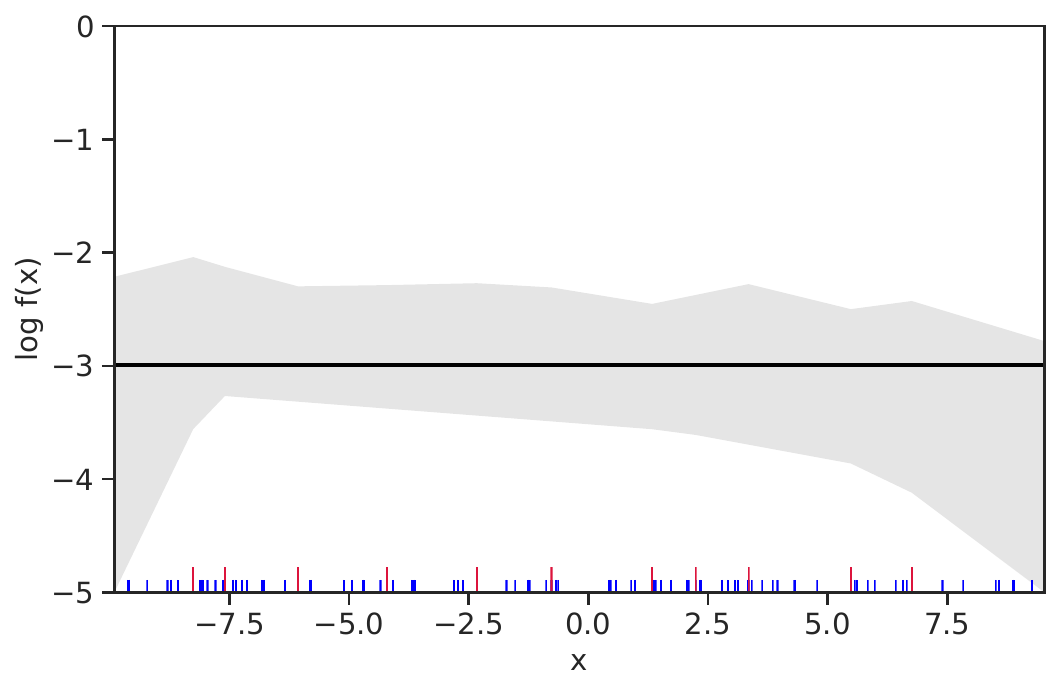} \hfill
\includegraphics[width=0.49\textwidth]{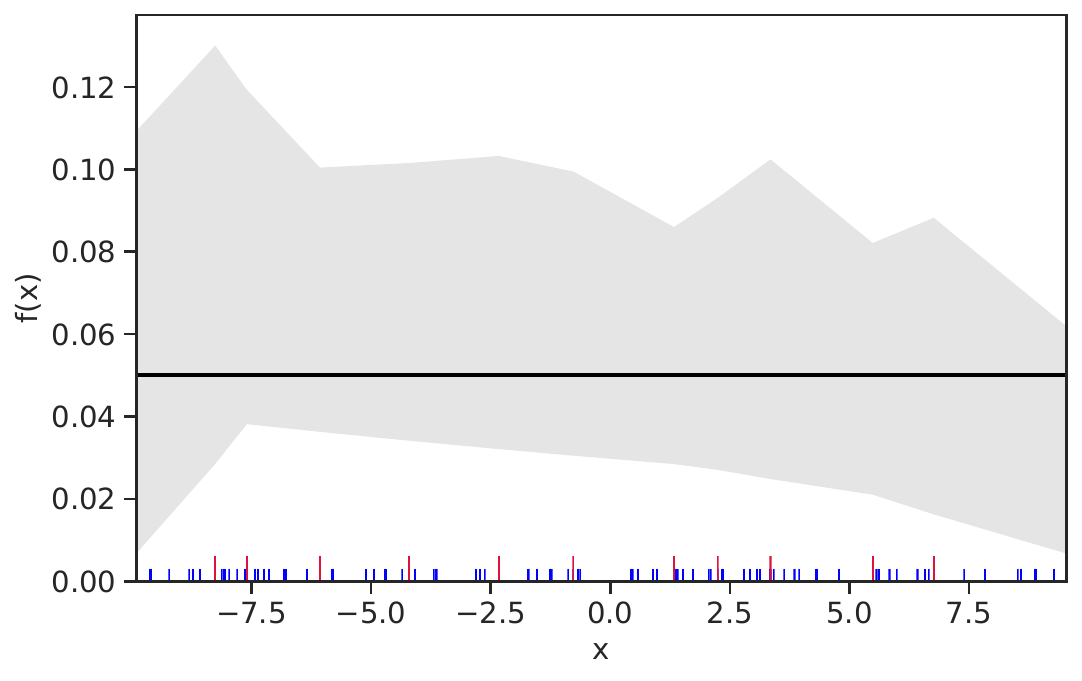} \\
\includegraphics[width=0.49\textwidth]{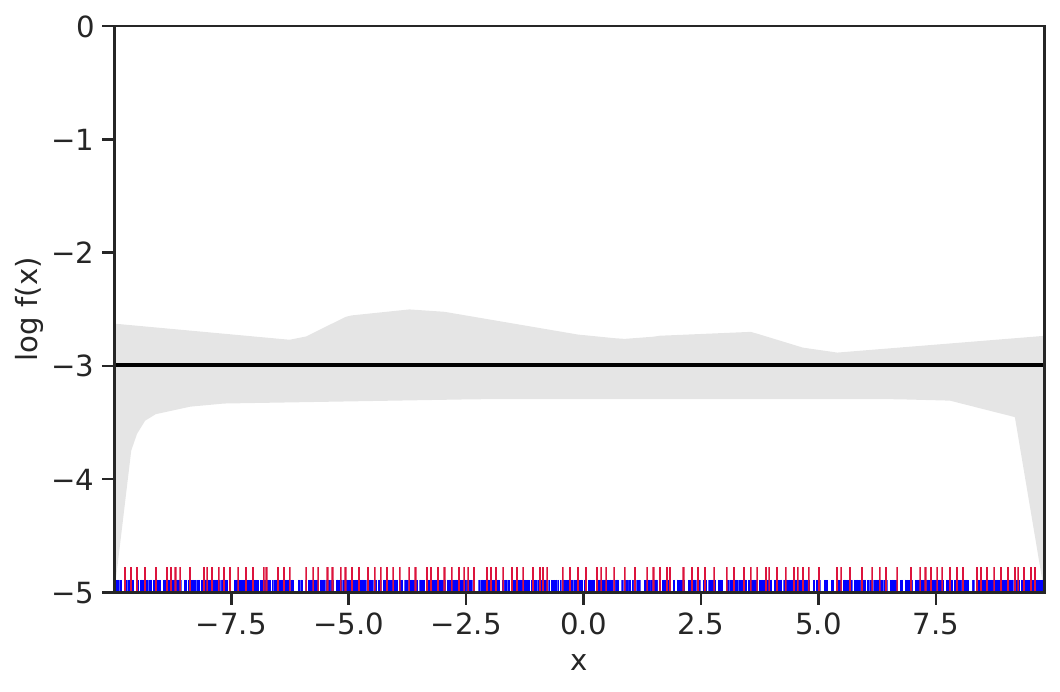} \hfill
\includegraphics[width=0.49\textwidth]{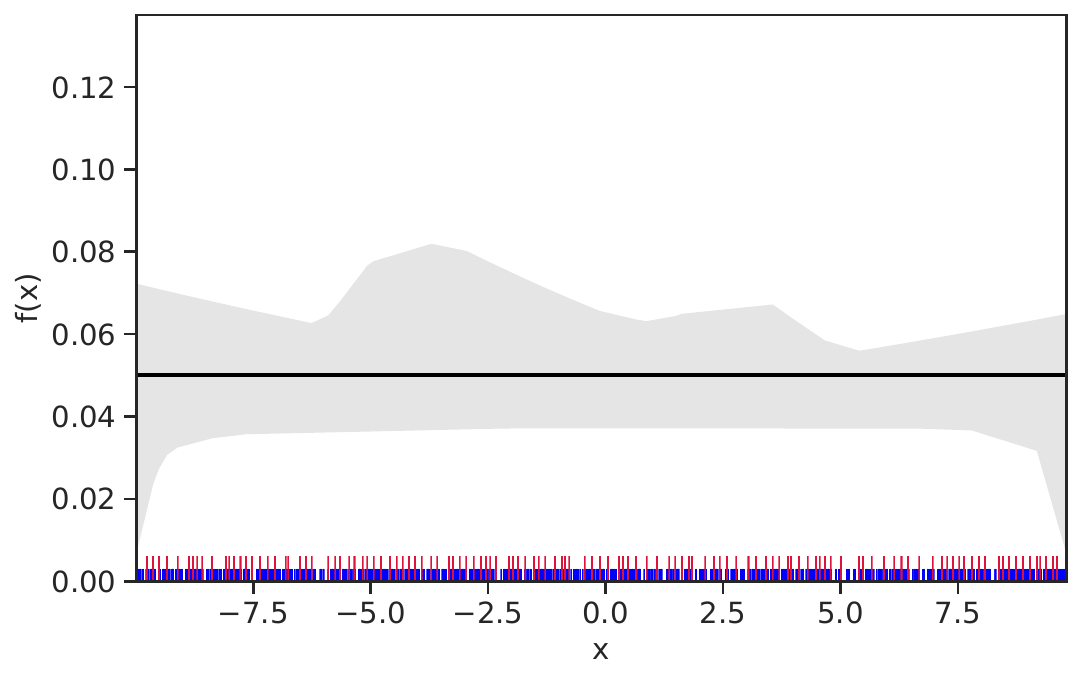} \\
\includegraphics[width=0.49\textwidth]{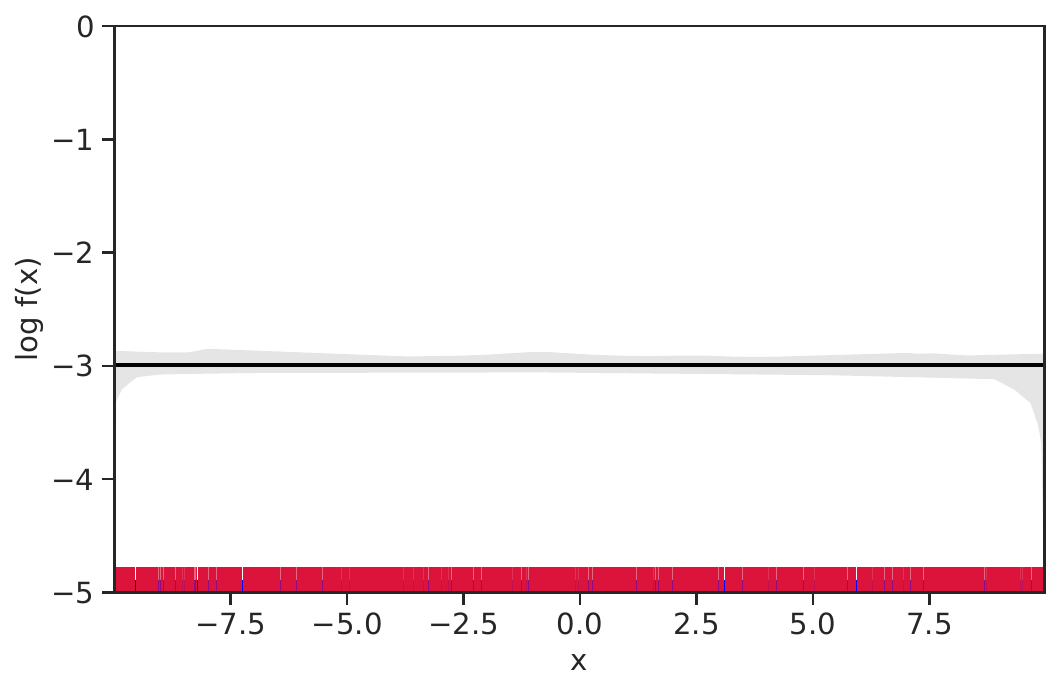} \hfill
\includegraphics[width=0.49\textwidth]{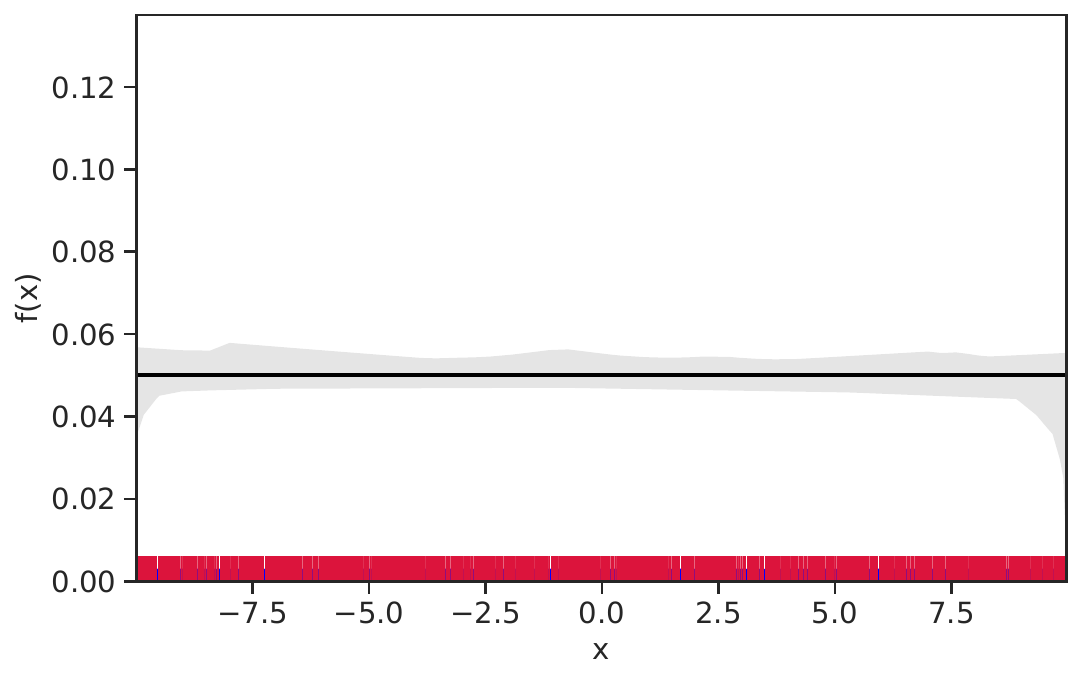}
% }
\vskip -0.2in
\caption{\textit{Confidence bands (shaded in gray) generated by Algorithm \ref{alg:ccp} and linear interpolation       
as described in Section~\ref{bands}, for a uniform density.  The left column shows the bands for the log density,   
while the right column shows the bands for the density. The solid black line marks the underlying (log)density.
Top, middle and bottom rows show results for sample sizes $n= 100$, $1000$, and $10000$, respectively.  
At the bottom of each plot, the observations
$X_i$ are indicated in blue (short lines), while the points $x_i, \; i=1,\ldots,m$ are marked in red (long lines).}}
\label{fig:bands:unif}
% \end{center}
\vskip -0.2in
\end{figure*}

\begin{figure*}[h!]
\vskip 0.2in
% \begin{center}
% \centerline{
\centering
\includegraphics[width=0.49\textwidth]{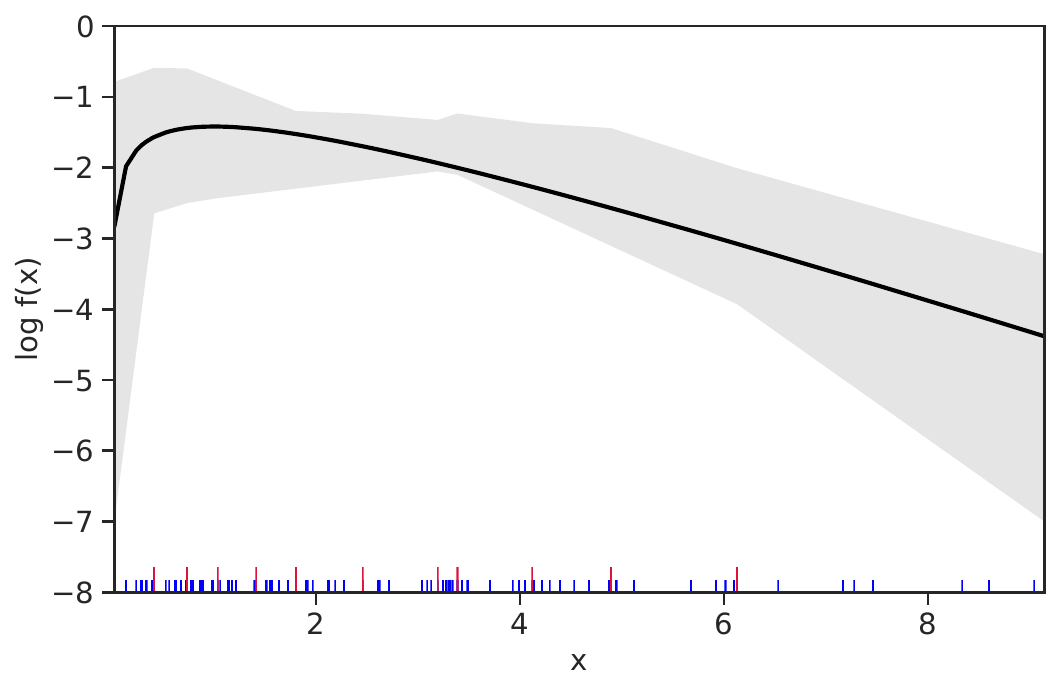} \hfill
\includegraphics[width=0.49\textwidth]{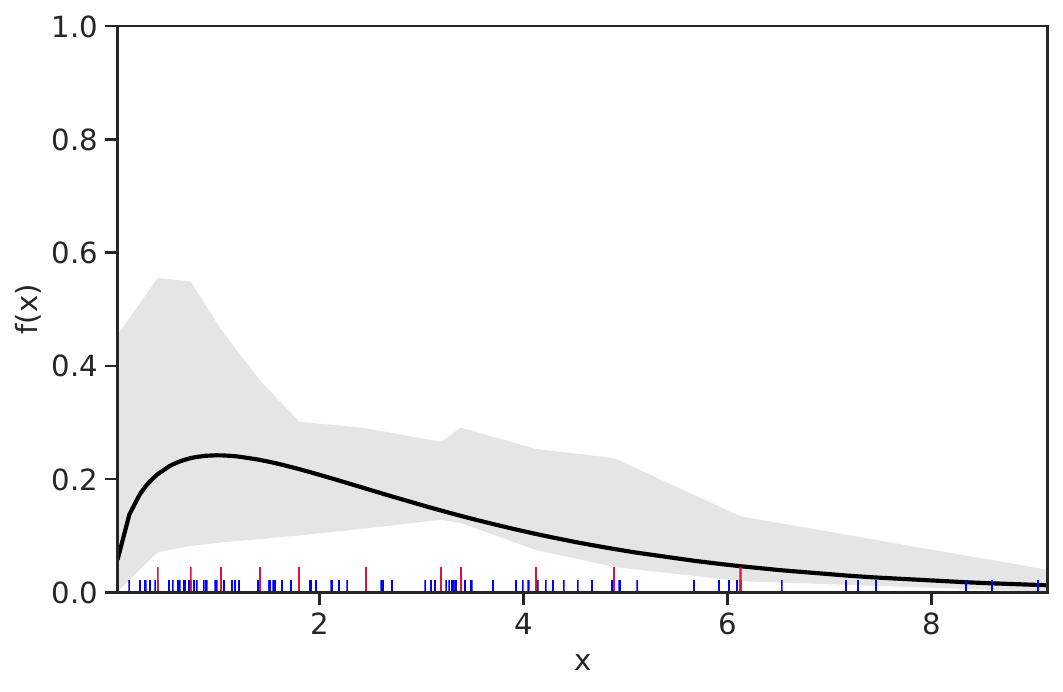} \\
\includegraphics[width=0.49\textwidth]{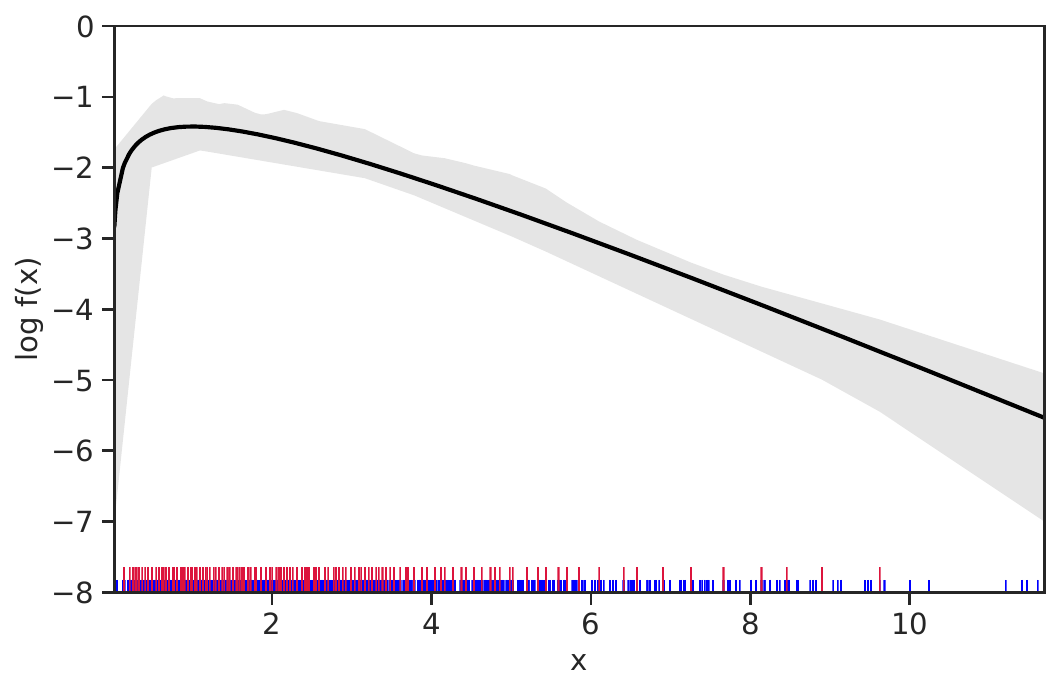} \hfill
\includegraphics[width=0.49\textwidth]{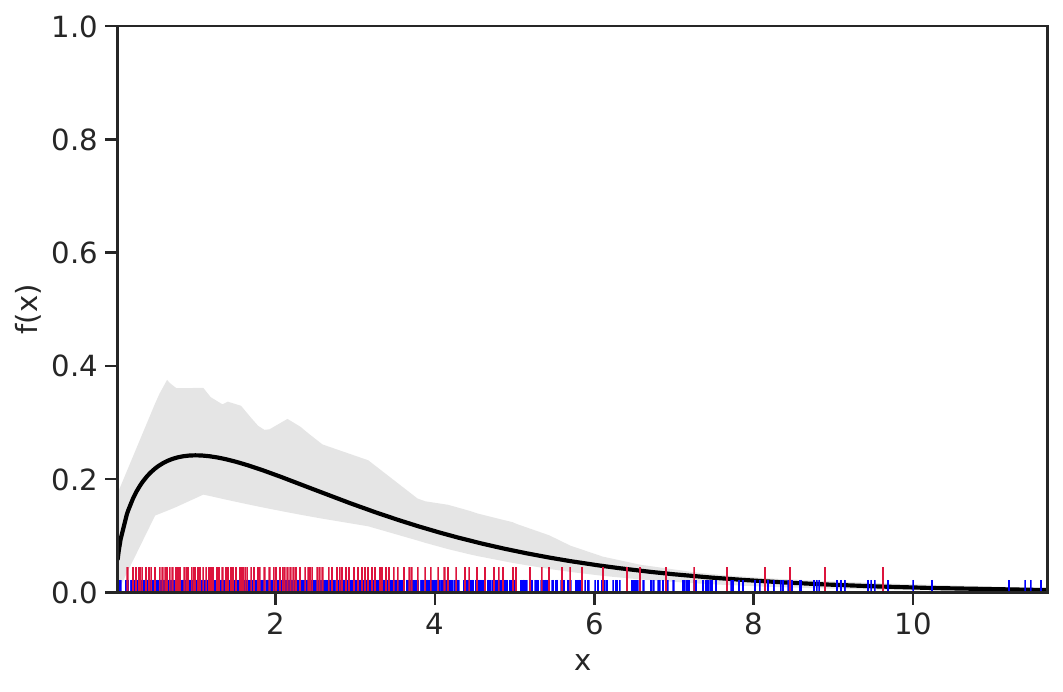} \\
\includegraphics[width=0.49\textwidth]{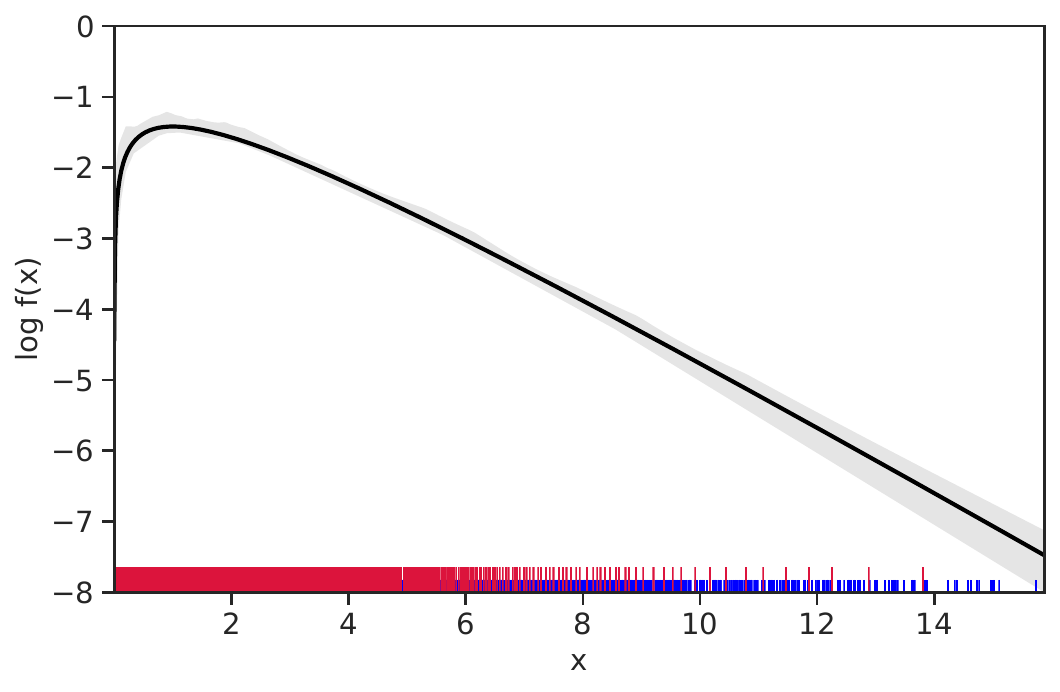} \hfill
\includegraphics[width=0.49\textwidth]{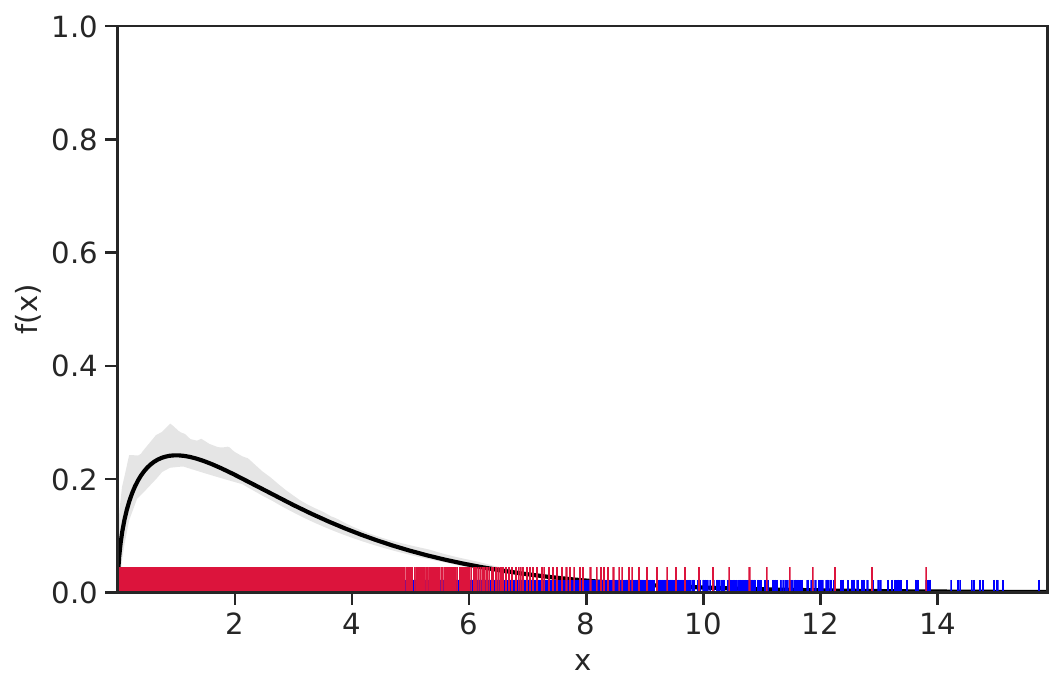}
% }
\vskip -0.2in
\caption{\textit{Confidence bands (shaded in gray) generated by Algorithm \ref{alg:ccp} and linear interpolation       
as described in Section~\ref{bands}, for a chi-squared density.  The left column shows the bands for the log density,   
while the right column shows the bands for the density. The solid black line marks the underlying (log)density.
Top, middle and bottom rows show results for sample sizes $n= 100$, $1000$, and $10000$, respectively.  
At the bottom of each plot, the observations
$X_i$ are indicated in blue (short lines), while the points $x_i, \; i=1,\ldots,m$ are marked in red (long lines).}}
\label{fig:bands:chisq}
% \end{center}
\vskip -0.2in
\end{figure*}

\begin{figure*}[h!]
\vskip 0.2in
% \begin{center}
% \centerline{
\centering
\includegraphics[width=0.49\textwidth]{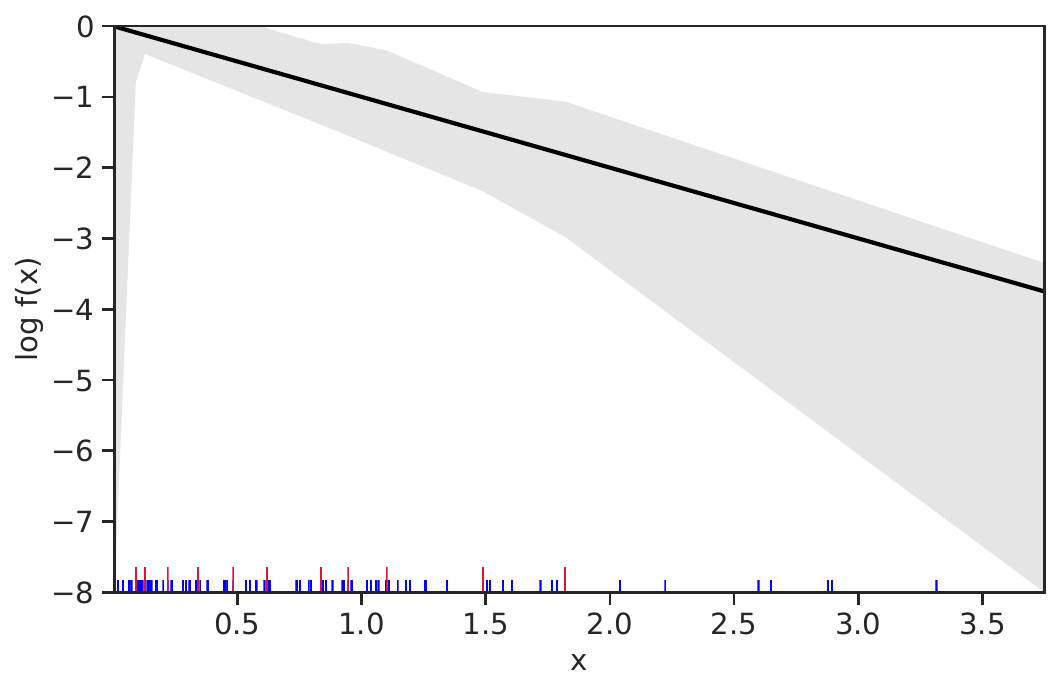} \hfill
\includegraphics[width=0.49\textwidth]{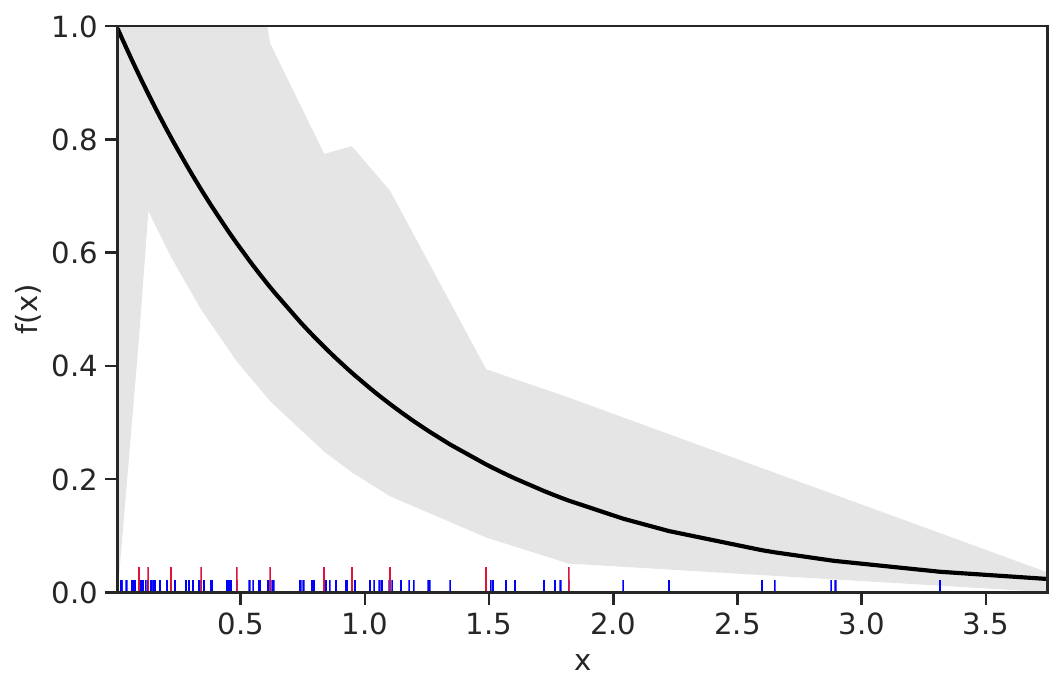} \\
\includegraphics[width=0.49\textwidth]{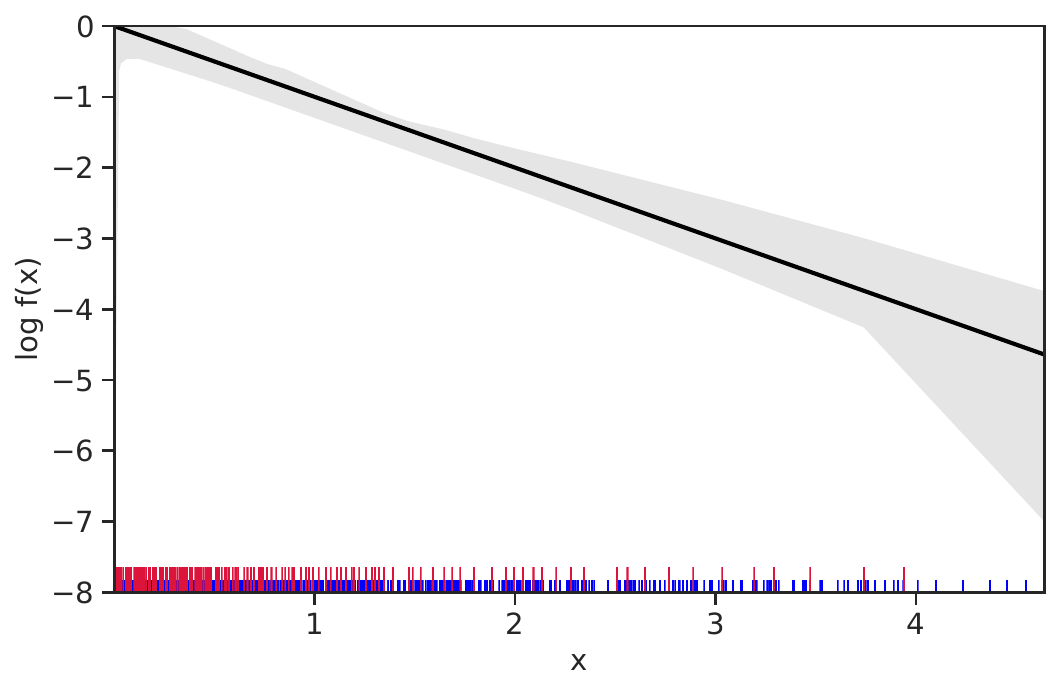} \hfill
\includegraphics[width=0.49\textwidth]{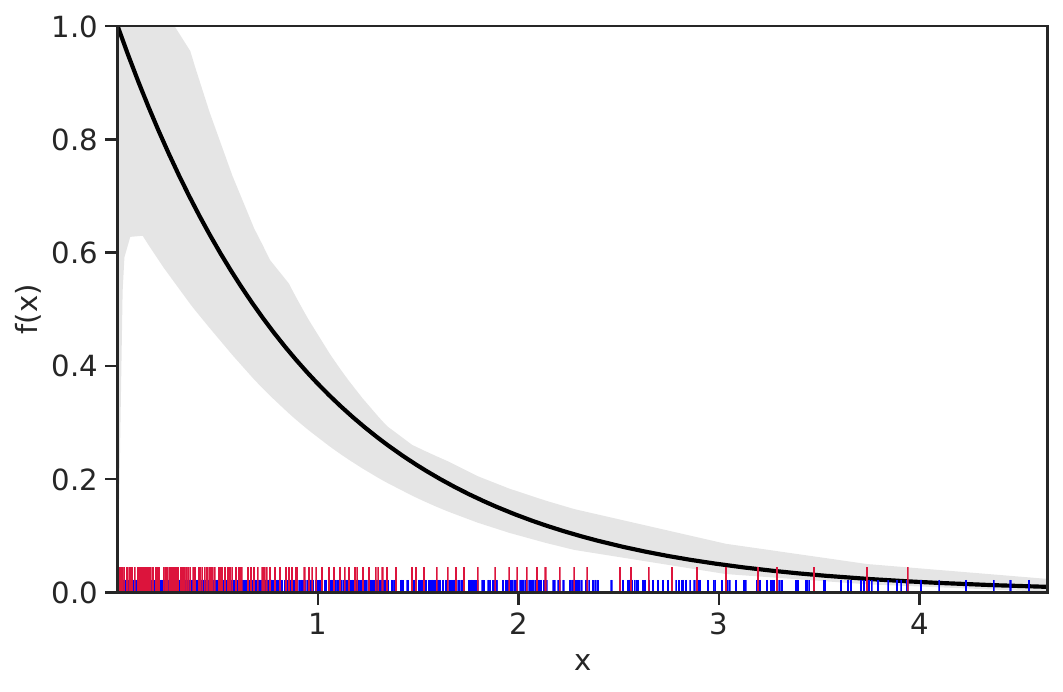} \\
\includegraphics[width=0.49\textwidth]{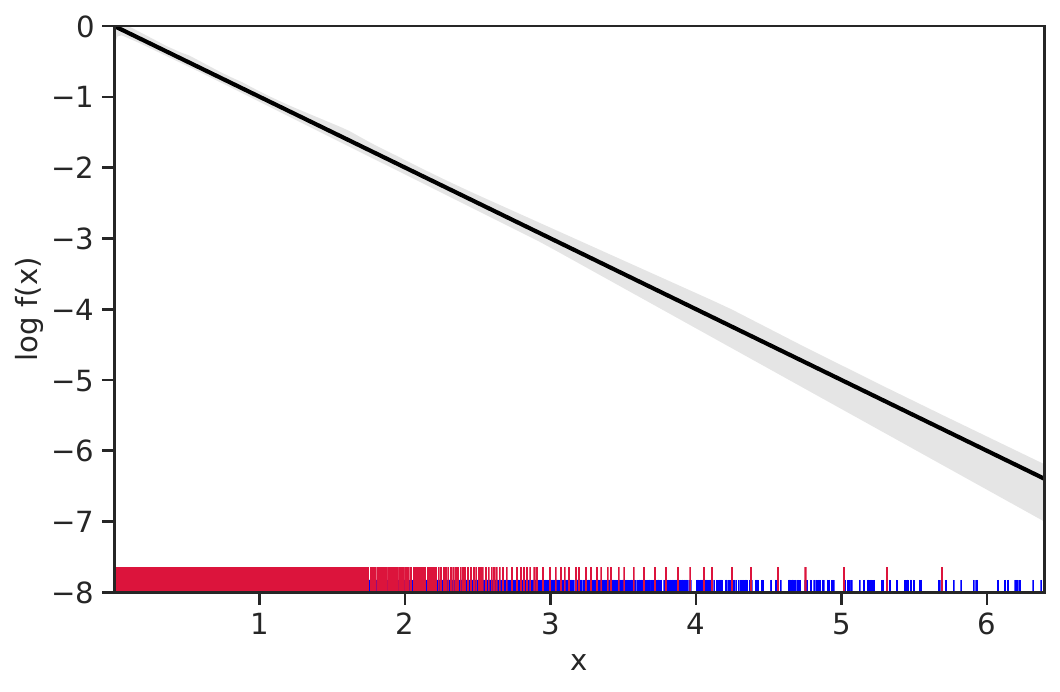} \hfill
\includegraphics[width=0.49\textwidth]{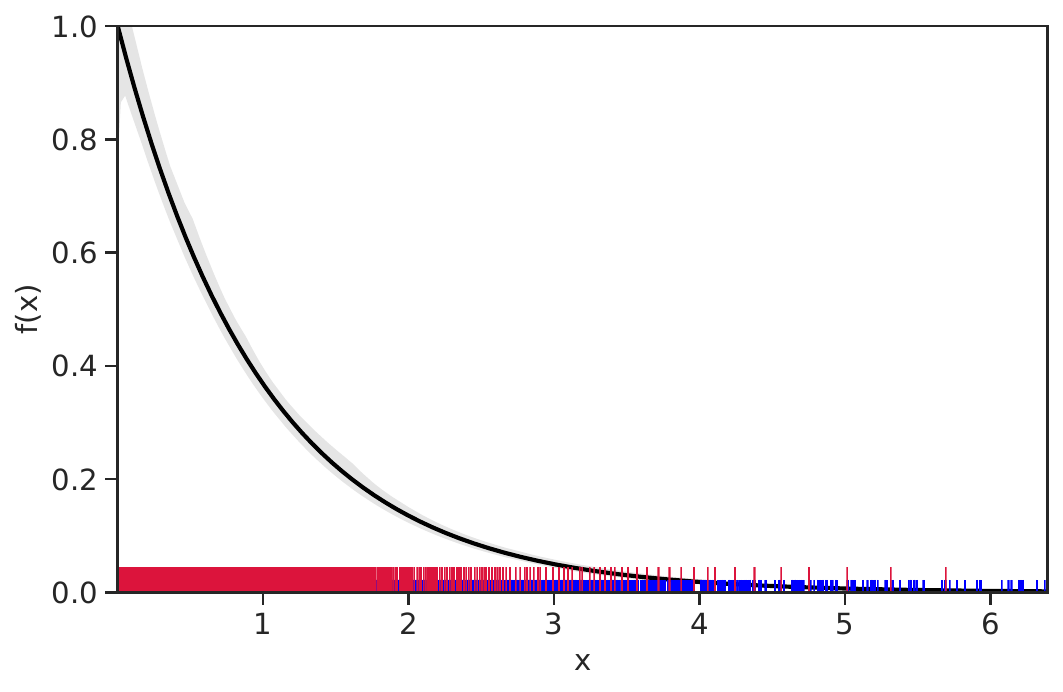}
% }
\vskip -0.2in
\caption{\textit{Confidence bands (shaded in gray) generated by Algorithm \ref{alg:ccp} and linear interpolation       
as described in Section~\ref{bands}, for an exponential density.  The left column shows the bands for the log density,   
while the right column shows the bands for the density. The solid black line marks the underlying (log)density.
Top, middle and bottom rows show results for sample sizes $n= 100$, $1000$, and $10000$, respectively.  
At the bottom of each plot, the observations
$X_i$ are indicated in blue (short lines), while the points $x_i, \; i=1,\ldots,m$ are marked in red (long lines).}}
\label{fig:bands:gam}
% \end{center}
\vskip -0.2in
\end{figure*}

\begin{figure*}[h!]
\vskip 0.2in
% \begin{center}
% \centerline{
\centering
\includegraphics[width=0.49\textwidth]{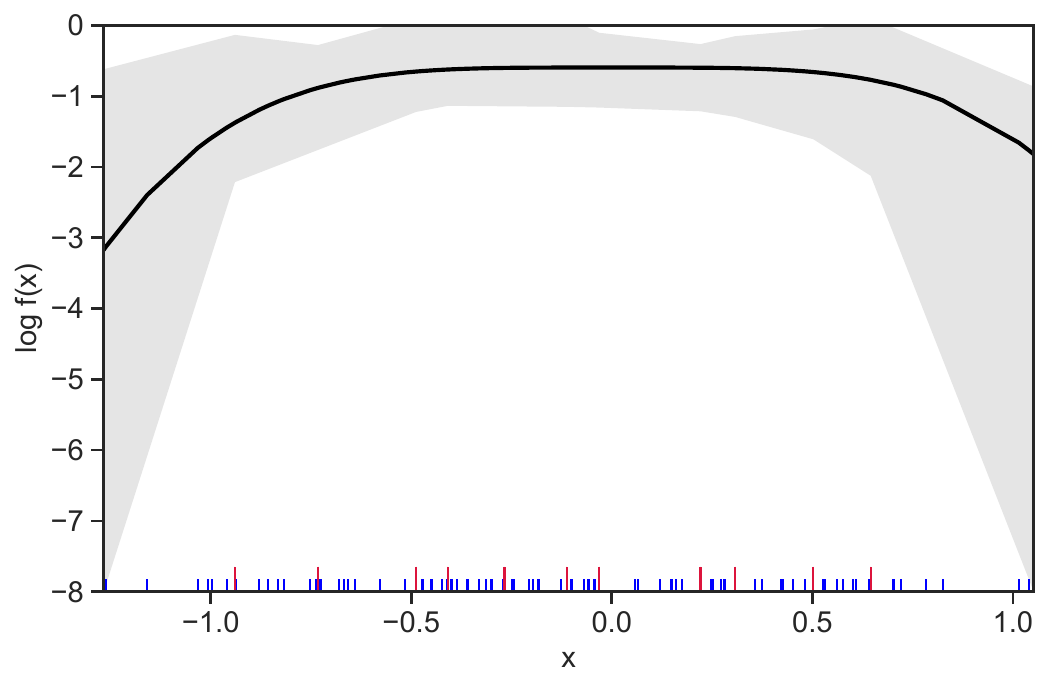} \hfill
\includegraphics[width=0.49\textwidth]{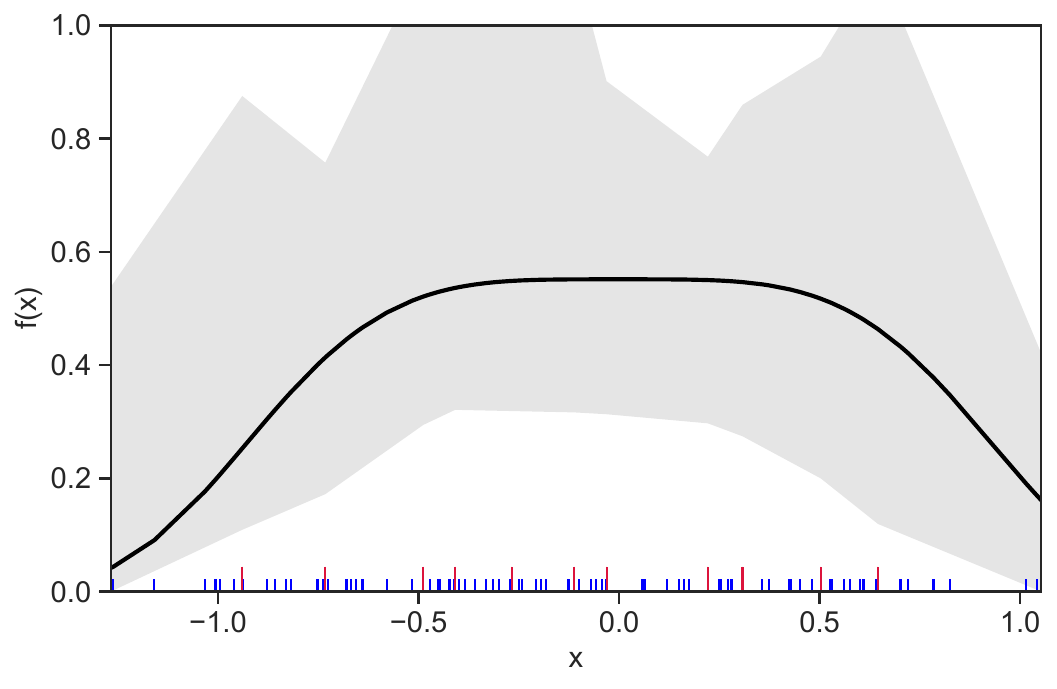} \\
\includegraphics[width=0.49\textwidth]{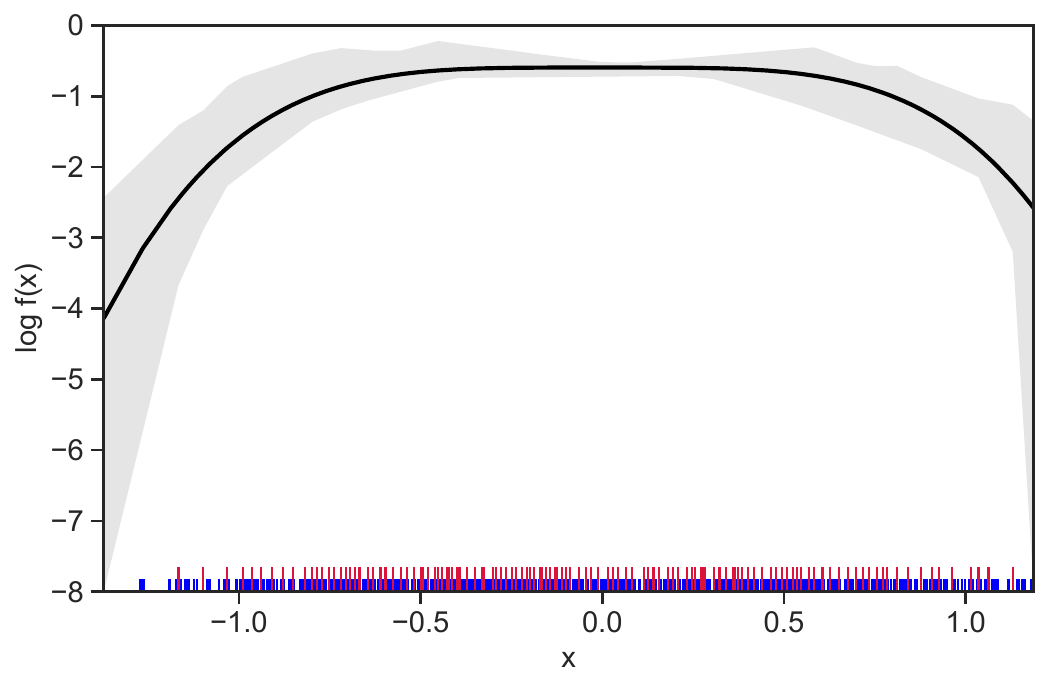} \hfill
\includegraphics[width=0.49\textwidth]{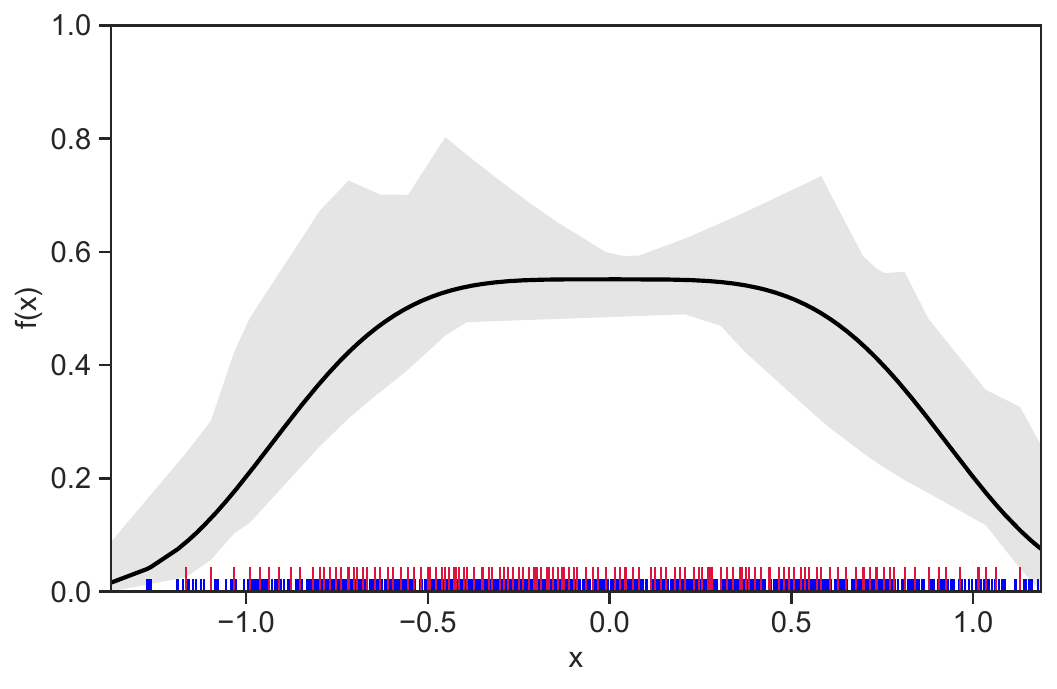} \\
\includegraphics[width=0.49\textwidth]{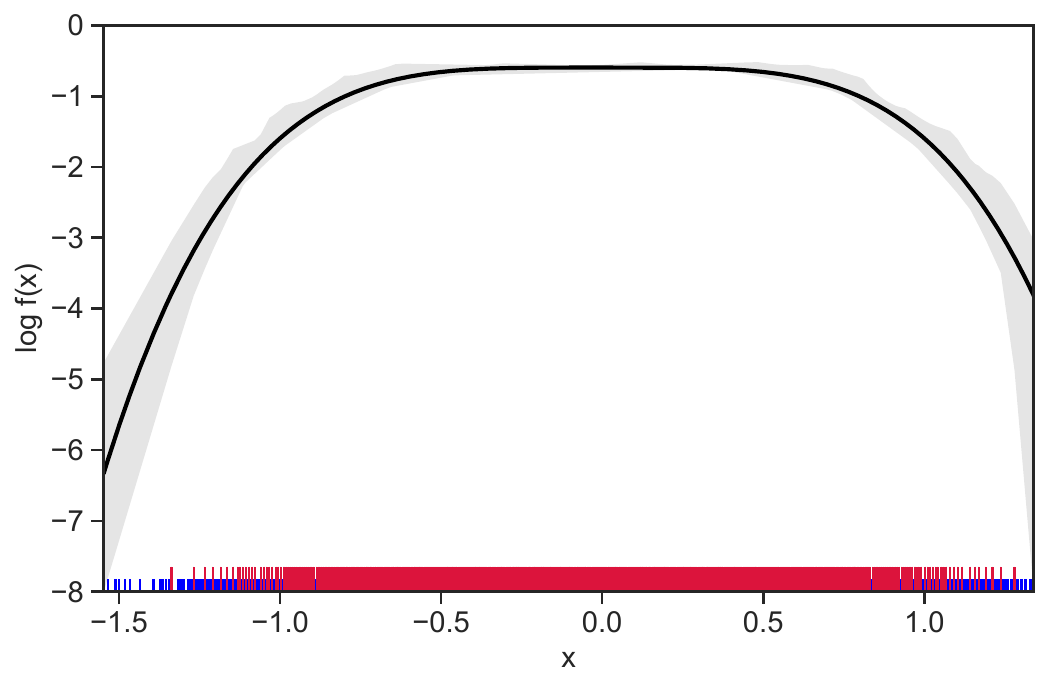} \hfill
\includegraphics[width=0.49\textwidth]{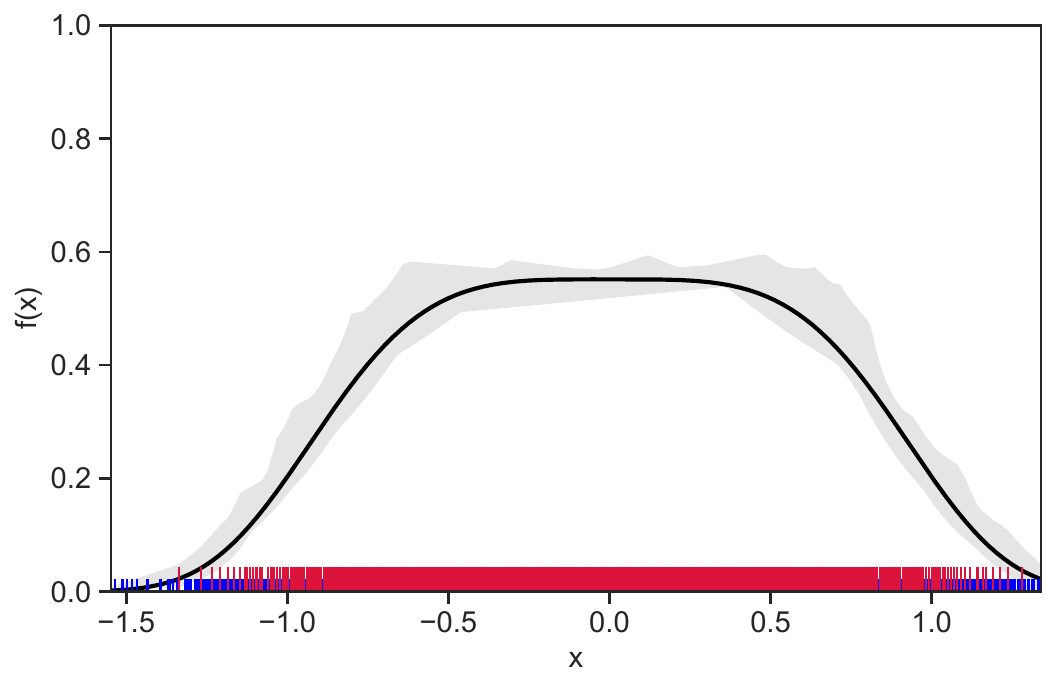}
% }
\vskip -0.2in
\caption{\textit{Confidence bands (shaded in gray) generated by Algorithm \ref{alg:ccp} and linear interpolation
as described in Section~\ref{bands}, for a density proportional to $\exp(-x^4)$.  
The left column shows the bands for the log density,
while the right column shows the bands for the density. The solid black line marks the underlying (log)density.
Top, middle and bottom rows show results for sample sizes $n= 100$, $1000$, and $10000$, respectively.
At the bottom of each plot, the observations
$X_i$ are indicated in blue (short lines), while the points $x_i, \; i=1,\ldots,m$ are marked in red (long lines).}}
\label{fig:bands:gennorm}
% \end{center}
\vskip -0.2in
\end{figure*}

Table \ref{tab:cvg_len_runtime} summarizes the results.  The table (reassuringly) shows us that the 
bands achieve coverage at or above the nominal level.  In Figures~\ref{fig:bands:gaus}--\ref{fig:bands:gennorm} 
we present a visualization of the bands, from a single repetition 
chosen at random, for each of the four underlying densities as well as for a density that is proportional
to $\exp (-x^4)$. In addition, we depict the bands for the case of
a larger sample with $n=10000$.  The figures and Table \ref{tab:cvg_len_runtime} show that while the bands are 
naturally wider when the sample size is small ($n=100$), they quickly tighten as the sample size grows 
($n \in \{1000,10000\}$).

As for the computational cost, we found that around 20-30 iterations of Algorithm \ref{alg:ccp} were enough to reach 
convergence, for each point $x_i, \; i=1,\ldots,m$. Table \ref{tab:cvg_len_runtime} shows that this translates into 
just a few seconds to compute the entire band when $n=100$, and a couple of minutes when $n=1000$.  We found that 
Algorithm \ref{alg:ccp} converged to the exact same solutions even when started from a number of different initial 
points, suggesting that it is in fact finding the global minimizers of the problems \eqref{Opt}.  Therefore, these 
runtimes appear to be reasonable, as it is worth bearing in mind that Algorithm \ref{alg:ccp} is effectively solving 
a potentially large number of non-convex optimization problems (precisely: 13, 39, and 188 such problems, corresponding 
to $n \in \{100,1000,10000\}$, respectively).  Moreover, we point out that the computation in Algorithm \ref{alg:ccp} 
can easily be parallelized, e.g. across the points $x_i, \; i=1,\ldots,m$.

In order to compare the confidence band with pointwise confidence intervals, 
we performed these experiments also with two  methods that compute pointwise 
confidence intervals for log-concave densities.  Azadbakhsh et al.~(2014) compare several such methods and report
that no one approach appears to uniformly dominate the others and that each method works well only in a certain range
of the data. The first group of methods examined by Azadbakhsh et al.~(2014) is based
on the pointwise asymptotic theory developed in Balabdaoui et al.~(2009) and requires
estimating a nuisance parameter, for which Azadbakhsh et al.~(2014)  investigate several options.
We picked the option that they report works best, namely their method (iv) in their Section~4. This method
is called `Asymptotic theory with approximation' in Table~\ref{tab:comparison}.
The second group of methods analyzed by Azadbakhsh et al.~(2014) concerns various  bootstrapping
schemes, and we chose the one
that they report to have the best performance, namely the ECDF-bootstrap, listed as (v) in their Section~4,
which we use with 250 bootstrap repetitions. This method computes the MLE for each bootstrap sample and then computes
the bootstrap percentile interval at a point $x_0$ based on the 250 bootstrap replicates of the MLE at $x_0$.
We used the R function \texttt{logConCI} produced by Azadbakhsh et al.~(2014) for implementing both methods.
With each of the two methods we compute the pointwise 90\% confidence interval for each point $x_0$ in the grid of points
that we use to evaluate empirical coverage. Since these are \textit{pointwise} 90\% confidence intervals, we expect
that the coverage for the band (i.e. the simultaneous coverage across all $x_0$ in the grid) is smaller than
90\%, but that the intervals are narrower than those for a simultaneous confidence band. This is confirmed
by the results in Table~\ref{tab:comparison}, which show that both methods seriously undercover.
The bootstrap is also seen
to be significantly more computationally intensive than the method based on asymptotic theory as well as 
Algorithm~\ref{alg:ccp}.

\begin{table}[h!]
%\vskip 0.15in
\begin{center}
\begin{tabular}{|l|c|c|ccc|c|ccc|r|}
\cline{3-11}
\multicolumn{2}{c|}{} & \multicolumn{4}{c|}{\textbf{Nominal coverage 90\%}} & 
\multicolumn{5}{c|}{\textbf{Nominal coverage 95\%}}\\
\hline
\multirow{2}{*}{\textbf{Distribution}} & \multirow{2}{*}{$n$} & \multirow{2}{*}{\textbf{Coverage}} & 
\multicolumn{3}{c|}{\textbf{Width}} & \multirow{2}{*}{\textbf{Coverage}} &
\multicolumn{3}{c|}{\textbf{Width}} &{\textbf{Time}} \\
% \abovespace\belowspace
& & & $Q_1$ & $Q_2$ & $Q_3$ & & $Q_1$ & $Q_2$ & $Q_3$ & (secs.)\\
\hline
%\abovespace
\multirow{2}{*}{Gaussian} & 100 & 0.96 & 0.58 & 0.65 & 0.61 & 0.98 & 0.66 & 0.74 & 0.70 &  5.6 \\
& 1000 & 0.95 & 0.24 & 0.28 & 0.24 & 0.97 & 0.27 & 0.31 & 0.27 & 151.2 \\
\hline
\multirow{2}{*}{Uniform} & 100 & 0.93 & 0.08 & 0.07 & 0.08 & 0.96 & 0.09 & 0.08 & 0.09 & 5.1 \\
& 1000 & 0.91 & 0.03 & 0.03 & 0.03 & 0.97 & 0.03 & 0.03 & 0.03 & 128.8 \\
\hline
\multirow{2}{*}{Chi-squared} & 100 & 0.94 & 0.36 & 0.28 & 0.19 & 0.98 & 0.41 & 0.33 & 0.23 & 5.4 \\
& 1000 & 0.94 & 0.16 & 0.12 & 0.07 & 0.97 & 0.18 & 0.13 & 0.08 & 132.5 \\
\hline
\multirow{2}{*}{Exponential} & 100 & 0.93 & 0.94 & 0.69 & 0.44 & 0.96 & 1.06 & 0.77 & 0.51 & 5.7 \\
& 1000 & 0.91 & 0.38 & 0.25 & 0.15 & 0.97 & 0.43 & 0.29 & 0.17 & 140.4 \\
\hline
\end{tabular} 
\end{center}
\vskip -0.1in
\caption{\textit{Empirical coverages, average widths at the sample quartiles (denoted $Q_1$, $Q_2$, and $Q_3$), 
and average runtimes in seconds, for the confidence bands generated by Algorithm \ref{alg:ccp} with linear interpolation
and nominal coverage levels 90\% and 95\%, 
when applied to $n \in \{100,1000\}$ observations drawn from four different underlying densities 
(Gaussian, uniform, chi-squared, and exponential).  All results are based on 1000 simulations.
The runtimes for the nominal coverage of 90\% are similar to the case of 95\% and are not displayed.}}
\label{tab:cvg_len_runtime}
\end{table}

\begin{table}[h!]
%\vskip 0.15in
\begin{center}
\begin{tabular}{|l|c|c|ccc|r|c|ccc|r|}
\cline{3-12}
\multicolumn{2}{c|}{} & \multicolumn{5}{c|}{\textbf{Asymptotic theory with approximation}} &
\multicolumn{5}{c|}{\textbf{Bootstrap}}\\
\hline
\multirow{2}{*}{\textbf{Distribution}} & \multirow{2}{*}{$n$} & \multirow{2}{*}{\textbf{Coverage}} &
\multicolumn{3}{c|}{\textbf{Width}} & {\textbf{Time}} & \multirow{2}{*}{\textbf{Coverage}} &
\multicolumn{3}{c|}{\textbf{Width}} &{\textbf{Time}} \\
% \abovespace\belowspace
& & & $Q_1$ & $Q_2$ & $Q_3$ & (secs.) & & $Q_1$ & $Q_2$ & $Q_3$ & (secs.)\\
\hline
%\abovespace
\multirow{2}{*}{Gaussian} & 100 & 0.22 & 0.13 & 0.20 & 0.13 & 1.4 & 0.37 & 0.16 & 0.16 & 0.12 & 98.2 \\
& 1000 & 0.08 & 0.04 & 0.08 & 0.04 & 6.4 & 0.15 & 0.04 & 0.07 & 0.03 & 194.8 \\
\hline
\multirow{2}{*}{Uniform} & 100 & 0.32 & 0.03 & 0.03 & 0.03 & 1.4 & 0.01 & 0.02 & 0.01 & 0.02 & 97.4 \\
& 1000 & 0.15 & 0.01 & 0.01 & 0.01 & 6.5 & 0.11 & 0.01 & 0.01 & 0.01 & 184.5 \\
\hline
\multirow{2}{*}{Chi-squared} & 100 & 0.06 & 0.08 & 0.04 & 0.02 & 1.4 & 0.52 & 0.04 & 0.03 & 0.02 & 97.2 \\
& 1000 & 0.00 & 0.02 & 0.01 & 0.00 & 6.6 & 0.31 & 0.02 & 0.01 & 0.01 & 191.0 \\
\hline
\multirow{2}{*}{Exponential} & 100 & 0.37 & 0.04 & 0.02 & 0.01 & 1.4 & 0.05 & 0.07 & 0.05 & 0.03 & 95.7 \\
& 1000 & 0.15 & 0.01 & 0.01 & 0.01 & 6.5 & 0.04 & 0.02 & 0.01 & 0.01 & 180.0 \\
\hline
\end{tabular}
\end{center}
\vskip -0.1in
\caption{\textit{Empirical coverages, average widths at the sample quartiles (denoted $Q_1$, $Q_2$, and $Q_3$),
and average runtimes in seconds, for the confidence bands obtained from pointwise asymptotic theory and approximation (iv)
in Azadbakhsh et al.~(2014) and by the bootstrap, with nominal level 90\%
in the same settings as in Table~\ref{tab:cvg_len_runtime}.}}
\label{tab:comparison}
\end{table}

% \clearpage
\clearpage

% \clearpage
\clearpage

% \clearpage
\clearpage

% \clearpage
\clearpage

% \clearpage
\clearpage

% [It's probably best to restrict plots of the bands to the range of the data $[x_1,x_m]$.]

\section{Discussion}

The paper shows how to construct confidence bands for a log-concave density by intersecting the log-concavity
constraint with an appropriate nonparametric confidence set. This approach has three strong points: First,
it produces confidence bands with a finite sample guaranteed
confidence level. Our simulations have shown that this guaranteed confidence level is not overly conservative.
Second, the approach allows to bring modern tools from optimization to bear on this problem. This aspect is
particularly important in a multivariate setting where it is known that computing the MLE is 
very time consuming. We expect that the key ideas of the univariate construction introduced here
can be carried over to the multivariate setting and we are working on implementing this program in the multivariate
setting. Third, it was shown that this approach results in attractive statistical properties, namely
that the confidence bands converge at nearly the parametric $n^{-\frac{1}{2}}$ rate when the log density
is $k$-affine. We conjecture that the width of these confidence bands will likewise achieve the optimal
minimax rate if $\log f$ is smooth rather than piecewise linear, and we leave the proof of such a result
as an open problem.

\section{Proofs}

\subsection{Proof of Lemma~\ref{integralbounds}}
It is well known that if $\log f(x)$ is a concave function then there exist real numbers $g_1,
\ldots, g_m$ such that
$$
\log f(x) \leq \ell_i +g_i \,(x-x_i)\ \ \mbox{ for } i=1,\ldots,m,
$$
see e.g. Boyd and Vandenberghe~(2004).
We found it computationally advantageous to employ these inequalities only for $i=2,\ldots,m-1$.
Those inequalities immediately yield
\be \label{conc2}
\ell_j \leq \ell_i +g_i \,(x_j-x_i)\ \ \mbox{ for all } i \in \{2,\ldots,m-1\},\ j\in \{i-1,i+1\},
\ee
as well as for $i=2,\ldots,m-1$:
$$
\int_{x_i}^{x_{i+1}} f(x)\,dx \leq \exp (\ell_i) \int_{x_i}^{x_{i+1}} \exp \bigl( g_i(x-x_i)\bigr)\,dx
= \exp (\ell_i) (x_{i+1}-x_i) \,E\bigl(g_i(x_{i+1}-x_i)\bigr)
$$
and  for $i=1,\ldots,m-2$:
$$
\int_{x_i}^{x_{i+1}} f(x)\,dx \leq \exp (\ell_{i+1}) \int_{x_i}^{x_{i+1}} \exp \bigl( g_{i+1}(x-x_{i+1})\bigr)\,dx
= \exp (\ell_{i+1}) (x_{i+1}-x_i) \,E \bigl( g_{i+1}(x_i-x_{i+1})\bigr).
$$
On the other hand, since $\log f$ is concave on $(x_i,x_{i+1})$ it cannot be smaller than the chord from 
$x_i$ to $x_{i+1}$. Hence for $i=1,\ldots,m-1$:
$$
\int_{x_i}^{x_{i+1}} f(x)\, dx \geq  \int_{x_i}^{x_{i+1}} \exp \Bl(\ell_i +(x-x_i)
\frac{\ell_{i+1} -\ell_i}{x_{i+1}-x_i} \Br) \,dx 
 =\ (x_{i+1}-x_i) \exp(\ell_i) \,E \bigl(\ell_{i+1}- \ell_i\bigr),
$$
proving the lemma. For later use we note the following fact:
\be \label{concext}
\mbox{If (\ref{conc2}) holds, then (\ref{conc2}) holds for all $j \in \{1,\ldots,m\}$, and the
$g_i$ are non-increasing in $i$.}
\ee
This is because (\ref{conc2}) implies both $\ell_{i+1} \leq \ell_i +g_i(x_{i+1}-x_i)$
and $\ell_i \leq \ell_{i+1}+g_{i+1}(x_i-x_{i+1})$, hence $g_i \geq \frac{\ell_{i+1}-\ell_i}{
x_{i+1}-x_i} \geq g_{i+1}$. This monotonicity property and (\ref{conc2}) yield for $j>i$:
$$
\ell_j \leq \sum_{k=i}^{j-1} g_k(x_{k+1}-x_k) +\ell_i \leq g_i \sum_{k=i}^{j-1} (x_{k+1}-x_k)+\ell_i
=\ell_i +g_i(x_j -x_i)
$$
and analogously for $j<i$. $\Box$

\subsection{Proof of Theorem~\ref{thm1}} Let $d>0$ be an integer that will be determined
later as a function of $f,I,J$, and set $B:=B_{max}-d$. So $d$ does not change with $n$
but $B$ does as $B_{max}$ increases with $n$. Define the event
$$
\AA_n(d):=\ \left\{ \left|  \int_{x_j}^{x_k} f(t)\,dt - p_{jk} \right| \ \leq \ \sqrt{\frac{p_{jk}
(1-p_{jk})}{n}} \sqrt{\log \log n} \ \mbox{ for all $(j,k) \in \II_B$} \right\}
$$
where $p_{jk}:= \frac{k-j}{n+1} 2^{s_n}$. We will prove the theorem with a sequence of lemmata.
The first lemma bounds the width of the confidence band 
on the discrete set $\{x_j, j \in \KK\}$, where $\KK :=\{j:\, j=1+i 2^B$ with $3 \leq i \leq n_B-3$
and $(x_{j-3\cdot 2^B},x_{j+3\cdot 2^B}) \subset J\}$:

\begin{Lemma} \label{B}
If $\{\ell_1,\ldots, \ell_m\,g_2,\ldots,g_{m-1}\}$ is a feasible point for the optimization
problem (\ref{Opt}), then on $\AA_n(d)$
$$
 \max_{j \in \KK} \left|\ell_j - \log f(x_j)\right| \ \leq \ 17
\sqrt{2^d \frac{\log \log n}{n}}.
$$
\end{Lemma}

In order to extend the bound over the discrete set $\{x_j,\ j \in \KK\}$ to a uniform bound
over the interval $I$ we will use the following fact, which is readily proved using elementary
calculations: 
\smallskip

If the linear function $L(t)$ and the concave function $g(t)$ satisfy $|L(t_i) -g(t_i)| \leq D$
for $t_1<t_2<t_3<t_4$, then $\sup_{t \in (t_2,t_3)} |L(t)-g(t)| \leq D \left(1 +2 \min \left(
\frac{t_3-t_2}{t_2-t_1}, \frac{t_3-t_2}{t_4-t_3} \right)\right)$.

\begin{Lemma} \label{C}
On the event $\AA_n(d)$
$$
\min \left( \max_{j \in \KK} \frac{x_{j+2^B}-x_j}{x_j -x_{j-2^B}},
\max_{j \in \KK} \frac{x_{j+2^B}-x_j}{x_{j+2\cdot 2^B} -x_{j+2^B}} \right) \leq \frac{9}{8}.
$$
\end{Lemma}

\nin Therefore any concave function $g$ for which $\{g(x_j), j\in \KK\}$ is feasible for (\ref{Opt})
must satisfy
$$
\sup_{x \in [x_{\underline{j}+2^B},x_{\overline{j}-2^B}]} \Bl|g(x) - \log f(x)\Br|\ \leq \ \frac{13}{4}
17 \sqrt{2^d \frac{\log \log n}{n}}
$$
on $\AA_n(d)$, where $\underline{j} := \min \KK$ and $\overline{j} := \max \KK$.
Hence this bound applies to the lower and upper confidence limits $\hat{\ell}$ and $\hat{\mu}$
since $\hat{\ell}$ is a concave function and $\{\hat{\ell}(x_i)$, $i=1,\ldots,m\}$ is feasible
for (\ref{Opt}), and for every real $x$ there exist a concave function $g$ such that $g(x)=\hat{\mu}(x)$
and $\{g(x_i)$, $i=1,\ldots,m\}$ is feasible for (\ref{Opt}). 

In order to conclude the proof
of the theorem we will show
\begin{Lemma} \label{A}
$$
\Pr \left\{\AA_n(d)\right\} \ \ra \ 1\ \mbox{ as $ n \ra \infty$}
$$
\end{Lemma}
and
\begin{Lemma} \label{D}
$$
\Pr \left\{ I \subset (x_{\underline{j}-2^B},x_{\overline{j}+2^B}) \right\} \ \ra \ 1\ \mbox{ as $ n \ra \infty$}
$$
for a certain $d=d(f,I,J)$.
\end{Lemma}
Then the claim of the theorem follows since the bound on the log densities
carries over to the densities as $\sup_{x \in I} |g(x) - \log f(x)| \leq 1$ implies
$$
\sup_{x \in I} \Bl|\exp (g(x)) -f(x)\Br|\ \leq \ 2 \Bl(\sup_{x \in I} f(x) \Br) \ \sup_{x \in I}
\Bl|g(x) -\log f(x)\Br|. 
$$ 
It remains to prove the Lemmata~\ref{B}--\ref{D}.
We will use the following two facts:
\begin{Fact} \label{fact1}
If $L(x)$ is a linear function, then $\int_{\alpha}^{\beta} \exp(L(x))\,dx= (\beta-\alpha)
\exp(L(\alpha))\,E\bigl(L(\beta) -L(\alpha)\bigr)$, where $E(\cdot)$ is given in Lemma~\ref{integralbounds}.
One readily checks that the function $\R^2 \ni (s,t) \mapsto \exp(t)\,E(s-t)$
is increasing in both $s$ and $t$. Furthermore, since $\frac{sinh(x)}{x}$ is increasing for $x>0$ 
we have for $s < t$ and $C>0$:
\begin{align*}
\exp(t+C)\,E\bigl(s-(t+C)\bigr) &= \frac{\exp(s) -\exp(t+C)}{s-(t+C)}\\
 &= \exp \left( \frac{s+t+C}{2}\right)
\frac{\sinh \left( \frac{t+C-s}{2}\right)}{\frac{t+C-s}{2}} \\
& \geq \exp \left( \frac{C}{2}\right) \,\exp \left( \frac{s+t}{2} \right)\,
\frac{\sinh \left( \frac{t-s}{2}\right)}{\frac{t-s}{2}} \\
& \geq \left( 1+ \frac{C}{2}\right) \frac{\exp(s) -\exp(t)}{s-t}\\
&= \left( 1+ \frac{C}{2}\right) \exp(t)\,E(s-t)
\end{align*}
and this inequality also holds for $s=t$ since $E$ is continuous.
\end{Fact}

\begin{Fact} \label{fact2}
Let $k\in (0,n+1)$ and $p:=\frac{k}{n+1}$. Then for $\alp \in (0,1)$
$$
\left.\begin{aligned}
q\textrm{Beta } (1-\alp, k, n+1-k) -p\\
-\Bl( q\textrm{Beta } (\alp, k, n+1-k) -p \Br)
\end{aligned}
\right\} \leq \sqrt{\frac{p(1-p)}{n+1}} \sqrt{2 \log \frac{1}{\alp}} + \frac{\log \frac{1}{\alp}}{n+1}
$$
\end{Fact}
where $q\textrm{Beta } (\alp, r,s)$ denotes the $\alp$-quantile of the beta distribution with
parameters $r$ and $s$. This fact follows from Propostion~2.1 in D\"{u}mbgen~(1998).
\medskip

\nin {\bf Proof of Lemma~\ref{B}:} Since $\log f$ is linear on $J$ we will assume that
it is non-increasing on $J$. The non-decreasing case is proved analogously.
Let $j\in \KK$ and set $k:=j+2^B$, so
$(j,k) \in \II_B$. The concavity constraint (\ref{CONC}) implies that the function
$h_{jk}(t):= \sum_{i=j}^{k-1} 1(x_i \leq t < x_{i+1}) \left( \ell_i +(t-x_i)\frac{
\ell_{i+1}-\ell_i}{x_{i+1}-x_i}\right)$ is concave on $(x_j,x_k)$ and hence not
smaller than its secant $g_{jk}(t):= \ell_j +(t-x_j)\frac{
\ell_{k}-\ell_j}{x_{k}-x_j}$. Hence
\begin{align}
\sum_{i=j}^{k-1} (x_{i+1}-x_i)  \exp(\ell_i)\,E(\ell_{i+1}-\ell_i) &=
\sum_{i=j}^{k-1} \int_{x_i}^{x_{i+1}} \exp \left(\ell_i +(t-x_i)\frac{
\ell_{i+1}-\ell_i}{x_{i+1}-x_i}\right) \, dt \nn \\
& = \int_{x_j}^{x_k} \exp \left( h_{jk}(t) \right) \, dt \nn \\
& \geq \int_{x_j}^{x_k} \exp \left( g_{jk}(t) \right) \, dt \nn \\
& = (x_k -x_j) \exp(\ell_j)\,E(\ell_{k}-\ell_j) \label{B1}
\end{align}

Suppose $\ell_j > \log f(x_j) +C$ where $C:=17 \cdot 2^{\frac{d}{2}} \sqrt{
\frac{\log \log n}{n}}$. If $\ell_k \geq \log f(x_k)$ then Fact~\ref{fact1} 
shows (note that $\log f(x_j) \geq \log f(x_k)$ as $f$ is non-increasing) that (\ref{B1}) 
is not smaller than
\begin{align}
&(x_k-x_j) \exp\bigl(\log f(x_j)+C\bigr)\,E\Bigl(\log f(x_k)-\left(\log f(x_j)+C\right)\Bigr)\nn \\
&\geq (x_k-x_j)  \left(1+\frac{C}{2}\right) \exp\left(\log f(x_j)\right)\,E\Bigl(\log f(x_k)-\log f(x_j)\Bigr) \nn\\
& = \left(1+\frac{C}{2}\right) \int_{x_j}^{x_k} f(t) \, dt \label{B2}
\end{align}
since $\log f(t)$ is linear on $(x_j,x_k)$. Now $p_{jk} \geq 2^{-d-4}$ since
$B=B_{max} -d$, and hence $C \geq \frac{17}{4} \sqrt{\frac{\log \log n}{p_{jk}n}}$.
So on the event $\AA_n(d)$, (\ref{B2}) is not smaller than
\begin{align*}
& \left( 1+ \frac{17}{8} \sqrt{\frac{\log \log n}{p_{jk}\,n}} \right) \left( p_{jk} -
\sqrt{\frac{p_{jk}(1-p_{jk})}{n}} \sqrt{\log \log n} \right) \\
& > p_{jk} +\frac{9}{8} \sqrt{\frac{p_{jk}(1-p_{jk})}{n}} \sqrt{\log \log n} \\
& > q\textrm{Beta }\left(1-\frac{\alp}{2(B+2)n_B t_n},(k-j)2^{s_n}, n+1-(k-j)2^{s_n} \right)\,=d_B
\end{align*}
by Fact~\ref{fact2}, since $B=B_{max}-d$ implies $n_B \leq \frac{n}{
2^{B_{max}+s_n-d}} \leq 16 \cdot 2^d$, while $t_n \leq \log (B_{max}+2) \leq \log \log_2 n$,
 and therefore $2(B+2) n_B t_n \leq
\alp (\log n)^2$ for $n$ large enough. Thus we arrive at a violation of the
constraint (\ref{UP}).

On the other hand, suppose $\ell_j > \log f(x_j) +C$ and $\ell_k < \log f(x_k)$.
Set $k=j+2^B$, $r:=j+2\cdot 2^B$ and $s:=j+3\cdot 2^B$. Then $x_r-x_k
\geq \frac{8}{9} (x_k-x_j)$ by (\ref{C2}). By (\ref{CONC}) and (\ref{concext})
$$
g_i \ \leq \ g_k \ \leq \ \frac{\ell_k -\ell_j}{x_k-x_j} \ \leq \ \frac{\log f(x_k) -\log f(x_j) -C}{
x_k-x_j}\ \ \ \mbox{ for $i=k,\ldots,s-1$}.
$$
Hence for $k \leq i <m\leq s-1$:
\be \label{C3}
g_i (x_m-x_i) \ \leq\ \log f(x_m) -\log f(x_i) - \frac{x_m-x_i}{x_k-x_j}C
\ee
since $\log f$ is linear on $(x_k,x_s)$. This yields
\begin{align*}
\ell_r & \leq \ell_k +g_k(x_r-x_k) \\
& \leq \log f(x_k) + \log f(x_r) -\log f(x_k) - \frac{x_r -x_k}{x_k-x_j}C \\
& \leq  \log f(x_r) - \frac{8}{9}C
\end{align*}
and likewise for $m=r+1,\ldots, s-1$:
\begin{align*}
\ell_m & \leq \ell_r +g_r (x_m-x_r) \\
& \leq \log f(x_r) -\frac{8}{9} C + \log f(x_m) -\log f(x_r) -\frac{x_m-x_r}{x_k-x_j}C \\
& \leq \log f(x_m) -\frac{8}{9} C
\end{align*}
Since the function $E(s)$ is positive and increasing for $ s\in \R$
we obtain with (\ref{C3}):
\begin{align*}
& \sum_{i=r}^{s-1} \exp(\ell_i) (x_{i+1}-x_i)\,E\left(g_i(x_{i+1}-x_i)\right) \\
& \leq \sum_{i=r}^{s-1}  \exp \left( \log f(x_i) -\frac{8}{9} C\right) 
 (x_{i+1}-x_i)\,E\Bigl(\log f(x_{i+1}) -\log f(x_i)\Bigr) \\
&= \exp \left(-\frac{8}{9}C\right)  \sum_{i=r}^{s-1} \int_{x_i}^{x_{i+1}} \exp \Bl(\log f(t) \Br) \,dt\\
&< \left(1-\frac{7}{9}C\right)  \int_{x_r}^{x_s} f(t)\,dt\ \ \ \ \ \mbox{ since 
$C\in \left[ \frac{17}{4} \sqrt{\frac{\log \log n}{p_{rs}n}}, \frac{1}{4}\right]$}\\
& \leq \left(1-\frac{7\cdot 17}{9\cdot 4} \sqrt{\frac{\log \log n}{p_{rs}n}}\right)
\left(p_{rs} + \sqrt{\frac{p_{rs}(1-p_{rs})}{n}} \sqrt{\log \log n} \right)\ \mbox{ on $\AA_n(d)$} \\
& < p_{rs} -2\sqrt{\frac{p_{rs}(1-p_{rs})}{n}} \sqrt{\log \log n}\\
& < q\textrm{Beta } \left( \frac{\alp}{2(B+2)n_B t_n}, (s-r)2^{s_n}, n+1-(s-r)2^{s_n} \right)\,=c_B
\end{align*}
by Fact~\ref{fact2}, yielding a violation of (\ref{DOWN2}). Therefore we conclude $\ell_j \leq
\log f(x_j) +C$ as claimed. The lower bound for $\ell_j$ follows analogously. $ \Box$
\medskip

\nin {\bf Proof of Lemma~\ref{C}:} $f$ is monotone on $J$ since it is log-linear. Consider first the case
where $f$ is non-increasing on $J$. On $\AA_n(d)$ we have for $j \in \KK$ and $k:=j+2^B$, $r:=j+2\cdot 2^B$:
%\be  \label{C1}
%\frac{\int_{x_{j+2^B}}^{x_{j+2\cdot 2^B}} f(t)\,dt}{\int_{x_j}^{x_{j+2^B}} f(t)\,dt} \geq
%\frac{p_{j+2^B,j+2\cdot 2^B} -\sqrt{\frac{p_{j+2^B,j+2\cdot 2^B} (1-p_{j+2^B,j+2\cdot 2^B})}{n}}
%\sqrt{\log \log n}}{p_{j,j+2^B}+\sqrt{\frac{p_{j,j+2^B} (1-p_{j,j+2^B})}{n}} \sqrt{\log \log n}} \geq \frac{8}{9}
%\ee
%for $n$ large enough, since $p_{j+2^B,j+2\cdot 2^B}=p_{j,j+2^B} \geq 2^{-d-4}$. Since $f$ is
%non-increasing on $J$, (\ref{C1}) implies
%\be \label{C2}
%x_{j+2^B}-x_j \leq \int_{x_j}^{x_{j+2^B}} \frac{f(t)}{f(x_{j+2^B})}\,dt \leq \frac{9}{8}
%\int_{x_{j+2^B}}^{x_{j+2\cdot 2^B}} \frac{f(t)}{f(x_{j+2^B})}\,dt \leq \frac{9}{8}
%(x_{j+2\cdot 2^B} -x_{j+2^B})
%\ee
\be  \label{C1}
\frac{\int_{x_{k}}^{x_{r}} f(t)\,dt}{\int_{x_j}^{x_{k}} f(t)\,dt} \geq
\frac{p_{kr} -\sqrt{\frac{p_{kr} (1-p_{kr})}{n}}
\sqrt{\log \log n}}{p_{jk}+\sqrt{\frac{p_{jk} (1-p_{jk})}{n}} \sqrt{\log \log n}} \geq \frac{8}{9}
\ee
for $n$ large enough, since $p_{kr}=p_{jk} \geq 2^{-d-4}$. Since $f$ is
non-increasing on $J$, (\ref{C1}) implies
\be \label{C2}
x_{k}-x_j \leq \int_{x_j}^{x_{k}} \frac{f(t)}{f(x_{k})}\,dt \leq \frac{9}{8}
\int_{x_{k}}^{x_{r}} \frac{f(t)}{f(x_{k})}\,dt \leq \frac{9}{8}
(x_{r} -x_{k})
\ee
and the claim of the Lemma follows. If $f$ is non-decreasing on $J$ then we obtain
analogously $x_{j+2^B}-x_j \leq \frac{9}{8} (x_j -x_{j-2^B})$ for all $j \in \KK$ and
the claim of the Lemma follows also. $\Box$
\medskip

\nin {\bf Proof of Lemma~\ref{A}:} Set $\alp_n :=(\log n)^{-\frac{1}{3}}$. Then for $(j,k) \in \II_B$
$$
\sqrt{\frac{p_{jk}(1-p_{jk})}{n+1}} \sqrt{2 \log \frac{1}{\alp_n}} + \frac{\log \frac{1}{\alp_n}}{n+1}
\ \leq \ \sqrt{\frac{p_{jk}(1-p_{jk})}{n}} \sqrt{\log \log n}
$$
for $n$ large enough as $p_{jk} =\frac{k-j}{n+1} 2^{s_n} =\frac{2^{B+s_n}}{n+1} \sim
\frac{2^{-d}}{8}$. There are $n_B \leq \frac{n}{2^{B_{max}+s_n-d}} \leq 16\cdot 2^d$ tuples
$(j,k)$ in $\II_B$ and
$$
\int_{x_j}^{x_k} f(t)\,dt \ \sim \ \textrm{ Beta }(2^{B+s_n},n+1-2^{B+s_n}).
$$
Hence, setting $j:=1$, $k:=1+2^B$ and using Fact~\ref{fact2}:
\begin{align*}
\Pr & \Bl\{ \AA_n(d)^c \Br\} \ \leq \ n_B \,\Pr \left\{ \left| \int_{x_j}^{x_k} f(t)\,dt -p_{jk} \right|
> \sqrt{\frac{p_{jk}(1-p_{jk})}{n}} \sqrt{\log \log n} \right\}\\
& \leq 16 \cdot 2^d \, \Pr \left\{ \int_{x_j}^{x_k} f(t)\,dt \not\in \Bl(q\textrm{Beta } ( \alp_n,
2^{B+s_n},n+1-2^{B+s_n}), q\textrm{Beta } ( 1-\alp_n,2^{B+s_n},n+1-2^{B+s_n})
\Br)\right\} \\
& =16 \cdot 2^d (2 \alp_n)\ \ra\ 0
\end{align*}
$\hfill \Box$
\medskip

\nin {\bf Proof of Lemma~\ref{D}:} Note that $J \setminus I =D_1 \cup D_2$, where $D_1$ is the
interval between the left endpoints of $I$ and $J$, and $D_2$ is the interval between the right
endpoints. Recall that $\underline{j}$ is the smallest index $j=1+i2^B$ such that
$x_{\underline{j}-3\cdot 2^B} \in J$. Hence $x_{\underline{j}-3\cdot 2^B}$ and
$x_{\underline{j}-2^B}$ must both fall into $D_1$ if at least $4 \cdot 2^B$ points $x_i$ fall
into $D_1$, i.e. if at least $4 \cdot 2^{B+s_n}$ observations $X_i$ fall into $D_1$.
Therefore
$$
\Pr \left\{ I \subset \left(x_{\underline{j}-2^B},x_{\underline{j}+2^B}\right) \right\}
\ \geq \ \Pr \left\{ F_n(D_i) \geq \frac{4 \cdot 2^{B+s_n}}{n},\ i=1,2\right\}
$$
where $F_n$ denotes the empirical distribution. Since $f$ is log-linear on $J$:
$$
q:= q(f,I,J):= \min \left( \int_{D_1} f(t)\, dt,\ \int_{D_2} f(t)\, dt \right)\ >\ 0.
$$
Therefore
$$
4 \frac{2^{B+s_n}}{n} -F(D_1) \ \leq \ 4 \frac{2^{\log_2 \frac{n}{8} -d}}{n} -q\ =\
2^{-1-d}-q\ <0
$$
for $d> \log_2 \frac{1}{2q}$, in which case Chebychev's inequality gives
$$
\Pr \left\{ F_n(D_1) <4 \frac{2^{B+s_n}}{n} \right\} \ \leq \ \frac{1}{4n(2^{-1-d}-q)^2}
\ \ra 0
$$
and the same result holds for $F_n(D_2)$. $\Box$
\medskip

\subsection*{References}

\begin{description}
%\item Axelrod, B., Diakonikolas, I., Sidiropoulos, A., Stewart, A. and Valiant, G. (2019).
%A polynomial time algorithm for log-concave maximum likelihood via locally exponential families.
%{\sl Adv. Neur. Inform. Process. Sys.}
\item Azadbakhsh, M., Jankowski, H. and Gao, X. (2014). Computing confidence intervals for
log-concave densities. {\sl Comput. Statist. Data Anal.} {\bf 75}, 248–264.
\item Balabdaoui, F., Rufibach, K. and Wellner, J. A. (2009). Limit distribution theory for maximum 
likelihood estimation of a log-concave density. {\sl Ann. Statist.} {\bf 37}, 1299–1331.
\item Birg\'{e}, L. (1997). Estimation of unimodal densities without smoothness assumptions. 
{\sl Ann. Statist.} {\bf 25}, 970–981.
\item Boyd, S. and Vandenberghe, L. (2004). {\sl Convex Optimization}. 
Cambridge University Press, New York.
\item Chan, H.P. and Walther, G. (2013). Detection with the scan and the average likelihood ratio.
{\sl Statist. Sinica} {\bf 23}, 409--428.
\item Cule, M., Samworth, R. and Stewart, M. (2010). Maximum likelihood
estimation of a multi-dimensional log-concave density. {\sl J. R. Stat. Soc. Ser.
B Stat. Methodol.} {\bf 72}, 545–607.
\item Davies, P.L. and Kovac, A. (2004). Densities, spectral densities and modality.
{\sl Ann. Statist.} {\bf 32}, 1093-1136.
\item Deng, H., Han, Q. and Sen, B. (2020). Inference for local parameters in convexity constrained models.
arXiv preprint arXiv:2006.10264.
\item Dinh, T. P. and Le Thi, H. A. (2014). Recent advances in DC programming and DCA. In {\sl 
Transactions on computational intelligence XIII}, 1-37. Springer, Berlin, Heidelberg.
\item Donoho, D.L. (1988). One-sided inference for functionals of a density.
{\sl Ann. Statist.} {\bf 16}, 1390-1420.
\item Doss, C. R. and Wellner, J. A. (2016). Global rates of convergence of the
MLEs of log-concave and s-concave densities. {\sl Ann. Statist.} {\bf 44}, 954–981.
\item D\"{u}mbgen, L. (1998). New goodness-of-fit tests and their application to nonparametric confidence sets.
{\sl Ann. Statist.} {\bf 26}, 288--314.
\item D\"{u}mbgen, L. (2003). Optimal confidence bands for shape-restricted curves.
{\sl Bernoulli} {\bf 9}, 423--449.
\item D\"{u}mbgen, L., H\"{u}sler, A. and Rufibach, K. (2007). Active set and EM algorithms
for log-concave densities based on complete and censored data. arXiv preprint arXiv:0707.4643v4.
\item D\"{u}mbgen, L. and Rufibach, K. (2009). Maximum likelihood estimation of a
log-concave density and its distribution function: basic properties and uniform
consistency. {\sl Bernoulli} {\bf 15}, 40–68.
\item Feng, O., Guntuboyina, A., Kim, A. K. H. and Samworth, R. J. (2020+). Adaptation in multivariate log-concave 
density estimation. {\sl Ann. Statist.}, to appear.
\item Hartman, P. (1959). On functions representable as a difference of convex functions. 
{\sl Pacific Journal of Mathematics} {\bf 9}(3), 707-713.
\item Hengartner, N.W. and Stark, P.B. (1995).  Finite-sample confidence envelopes for shaperestricted densities. 
{\sl Ann. Statist.} {\bf 23}, 525–550.
\item Horst, R., Pardalos, P. M. and Van Thoai, N. (2000). {\sl Introduction to global optimization.}
 Springer Science and Business Media.
\item Horst, R. and Thoai, N. V. (1999). DC programming: overview. 
{\sl Journal of Optimization Theory and Applications} {\bf 103}(1), 1-43.
\item Khamaru, K. and Wainwright, M. J. (2018). Convergence guarantees for a class of non-convex and non-smooth 
optimization problems. arXiv preprint arXiv:1804.09629.
\item Kim, A. K. H. and Samworth, R. J. (2016). Global rates of convergence in
log-concave density estimation. {\sl Ann. Statist.} {\bf 44}, 2756–2779.
\item Kim, A. K. H., Guntuboyina, A. and Samworth, R. J. (2018). Adaptation in log-concave density
estimation. {\sl Ann. Statist.} {\bf 46}, 2279–2306.
\item Le Thi, H. A. and Dinh, T. P. (2014). DC programming and DCA for general DC programs. In {\sl 
Advanced Computational Methods for Knowledge Engineering}, 15-35. Springer, Cham.
\item Li, H., Munk, A., Sieling, H. and Walther, G. (2020). The essential histogram.
{\sl Biometrika} {\bf 107}, 347–364. 
\item Lipp, T. and Boyd, S. (2016). Variations and extension of the convex–concave procedure. 
{\sl Optimization and Engineering} {\bf 17}(2), 263-287.
\item Liu, Y. and Wang, Y. (2018). A  fast  algorithm  for  univariate  log-concave  densityi estimation.
{\sl  Aust. N. Z. J. Stat.} {\bf 60}(2), 258–275.
\item Pal, J. K., Woodroofe, M. and Meyer, M. (2007). Estimating a Polya frequency function.
 In {\sl Complex datasets and inverse problems. IMS Lecture Notes Monogr. Ser.} {\bf 54}, 239–249. 
Inst. Math. Statist., Beachwood, OH.
%\item Rathke, F. and Schn\"{o}rr, C. (2019). Fast multivariate log-concave density estimation.
%{\sl Comput. Statist. Data Anal.} {\bf 140}, 41-58.
\item Rivera, C. and Walther, G. (2013). Optimal detection of a jump in the intensity of a Poisson
process or in a density with likelihood ratio statistics. {\sl Scand. J. Stat.} {\bf 40}, 752–769.
\item Samworth, R.J. (2018). Recent progress in log-concave density estimation. {\sl Statist. Sci.} {\bf 33}, 493-509.
\item Saumard, A. and Wellner, J.A. (2014).  
Log-concavity and strong log-concavity: a review.
{\sl Statistics Surveys} {\bf 8}, 45-114.
\item Schuhmacher, D. and D\"{u}mbgen, L. (2010). Consistency of multivariate log-concave density estimators. 
{\sl Statist. Probab. Lett.} {\bf 80}, 376–380. 
\item Seregin, A. and Wellner, J. A. (2010). Nonparametric estimation of multivariate convex-transformed densities. 
{\sl Ann. Statist.} {\bf 38}, 3751–3781.
\item Shorack, G.R. and Wellner, J.A. (1986). {\sl Empirical Processes with Applications to Statistics}. 
Wiley, New York.
\item Smola, A. J., Vishwanathan, S. V. N. and Hofmann, T. (2005). Kernel Methods for Missing Variables. 
In {\sl AISTATS}.
\item Sriperumbudur, B. K. and Lanckriet, G. R. (2009, December). On the convergence of the concave-convex procedure. 
In {\sl Proceedings of the 22nd International Conference on Neural Information Processing Systems},  1759-1767. 
Curran Associates Inc.
\item Tao, P. D. (1986). Algorithms for solving a class of nonconvex optimization problems. Methods of subgradients. 
In {\sl North-Holland Mathematics Studies}  {\bf 129},  249-271. North-Holland.
\item Walther, G. (2002). Detecting the presence of mixing with multiscale maximum likelihood. 
{\sl J. Amer. Statist. Assoc.} {\bf 97}, 508–513.
\item Walther, G. (2009). Inference and modeling with log-concave distributions.
{\sl Statist. Sci.} {\bf 24}, 319–327.
\item Walther, G. (2010). Optimal and fast detection of spatial clusters with scan statistics.
{\sl Ann. Statist.} {\bf 38}, 1010–1033.
\item Walther, G. and Perry, A. (2019). Calibrating the scan statistic: finite sample performance vs. asymptotics.
arXiv preprint arXiv:2008.06136.
\item Yuille, A. L. and Rangarajan, A. (2003). The concave-convex procedure. 
{\sl Neural computation} {\bf 15}(4), 915-936.
\end{description}

\end{document}